\documentclass[namedate,imamat]{ima-authoring-template}

\usepackage{bm}
\usepackage{float}
\usepackage{amsmath}
\usepackage[authoryear]{natbib}

\theoremstyle{remark}
\newtheorem{remark}{Remark}
\numberwithin{equation}{section}

\newcommand{\supa}{^{\scriptscriptstyle (\alpha)} \!}
\newcommand{\supap}{^{\scriptscriptstyle (\alpha) \prime} \!}
\newcommand{\supo}{^{\scriptscriptstyle \! (1)} \!}

\begin{document}

\DOI{DOI HERE}
\copyrightyear{2026}
\vol{00}
\pubyear{2026}
\access{Advance Access Publication Date: Day Month Year}
\appnotes{Paper}
\firstpage{1}

\title[Construction of Laguerre pseudospectral differentiation matrices]{Construction of Laguerre pseudospectral differentiation matrices}
\author{Emma Nel*~\ORCID{0000-0003-3998-8899}
\address{\orgdiv{Department of Mathematical Sciences}, \orgname{Stellenbosch University}, \orgaddress{\street{Stellenbosch}, \postcode{7600}, \country{South Africa}}}}
\author{Nicholas Hale~\ORCID{0000-0002-2023-0044}
\address{\orgdiv{Department of Mathematical Sciences}, \orgname{Stellenbosch University}, \orgaddress{\street{Stellenbosch}, \postcode{7600}, \country{South Africa}}}}

\authormark{E. Nel and N. Hale}

\corresp[*]{Corresponding author: \href{email:email-id.com}{eanel@sun.ac.za}}

\received{Date}{0}{Year}
\revised{Date}{0}{Year}
\accepted{Date}{0}{Year}

\abstract{In this paper, we present a stable and efficient approach for constructing Laguerre pseudospectral differentiation matrices. The proposed method reformulates the off-diagonal entries and computes all required quantities simultaneously using an existing fast algorithm that also generates the collocation nodes. For the diagonal entries, a closed-form expression is employed to improve numerical accuracy. This construction avoids the catastrophic cancellation present in classical formulations and yields an all-in-one procedure for generating differentiation matrices. Numerical experiments demonstrate improved robustness and sustained high accuracy for significantly larger numbers of collocation points compared to standard implementations.}


\maketitle

\section{Introduction}
\label{sec:intro}

Spectral and pseudospectral methods based on classical orthogonal polynomials have a long history of efficiency, robustness, and high accuracy in the numerical solution of differential equations~\citep{boyd2001,canuto1988,gottlieb1977}. When such methods are extended to semi-infinite domains, the Laguerre polynomials arise as a natural choice of basis, owing to their orthogonality on $[0,\infty)$ with respect to an exponentially decaying weight. As a result, Laguerre-based spectral and pseudospectral methods have been widely used for ordinary and partial differential equations posed on the half-line~\citep{huang2024,wang2008,wang2016,gheorghiu2013}.

The use of Laguerre polynomials in spectral methods dates back to the 1970s, beginning with the work of~\citet{gottlieb1977}, and was further developed in the 1980s by~\citet{maday1985} and~\citet{canuto1988}. Since then, Laguerre spectral and pseudospectral methods have been applied to a variety of problems on semi-infinite domains, including fluid dynamics, heat transfer, quantum mechanics, and boundary value problems \citep{wang2008,wang2016,gheorghiu2013}. Theoretical foundations of these methods are detailed in textbooks by~\citet{boyd2001},~\citet{canuto2006}, and~\citet{shen2011}. Despite this long history, Laguerre spectral methods remain comparatively less mature than spectral methods on bounded domains, such as Chebyshev and Legendre methods, and Fourier methods on periodic domains. Several aspects of their theoretical and numerical behaviour continue to attract active research, including stability analysis, numerical robustness, and the optimal choice of scaling parameters \citep{huang2024,wang2024}.

A well-known practical difficulty in Laguerre-based methods is the rapid growth of high-degree Laguerre polynomials and their derivatives, which can lead to severe numerical instability in floating-point arithmetic~\citep{funaro1990,shen2000,huang2024}. In particular, classical implementations may suffer from overflow or underflow as the polynomial degree increases, resulting in a loss of accuracy or even complete failure of the computation. These issues affect not only the evaluation of basis functions themselves, but also the numerical procedures built upon them, such as root-finding algorithms, quadrature rules, barycentric interpolation formulas, and the construction of differentiation matrices.

Among these applications, the construction of pseudospectral differentiation matrices is particularly important in collocation methods. Although the differential operators they approximate are well-behaved, standard formulas for constructing Laguerre pseudospectral differentiation matrices often involve intermediate quantities whose magnitudes exceed the limits of double-precision arithmetic, leading to overflow. Consequently, numerical instability may arise for a relatively modest number of collocation points, limiting the practical range of Laguerre pseudospectral methods despite their favourable approximation properties. 

\begin{figure}[h]
    \centering
    \includegraphics[width=0.48\linewidth]{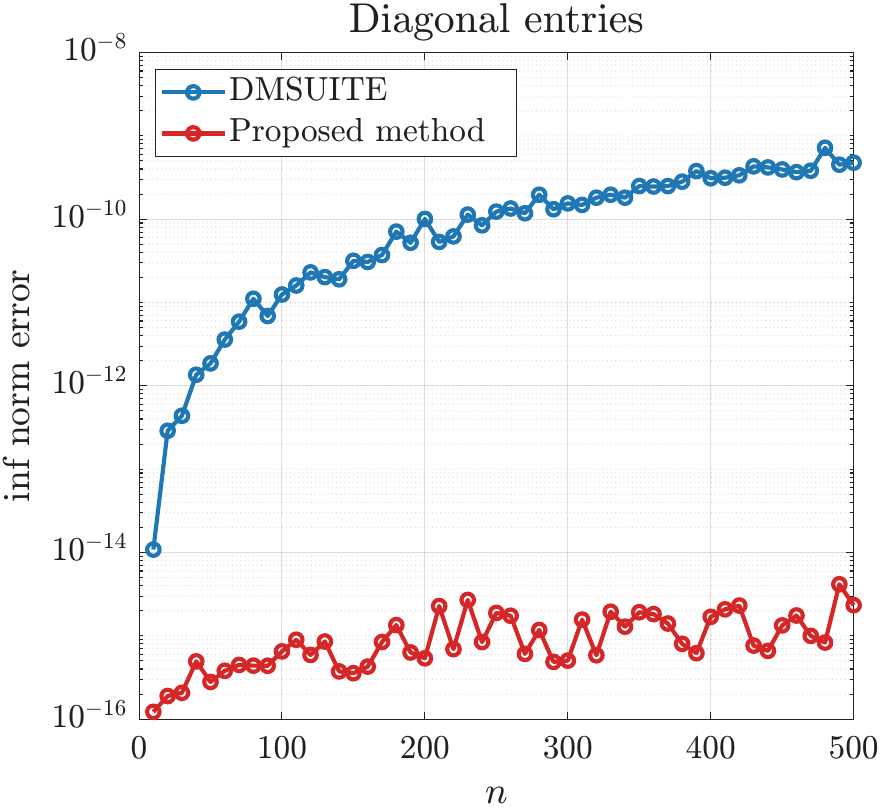}
    \hfill
    \includegraphics[width=0.48\linewidth]{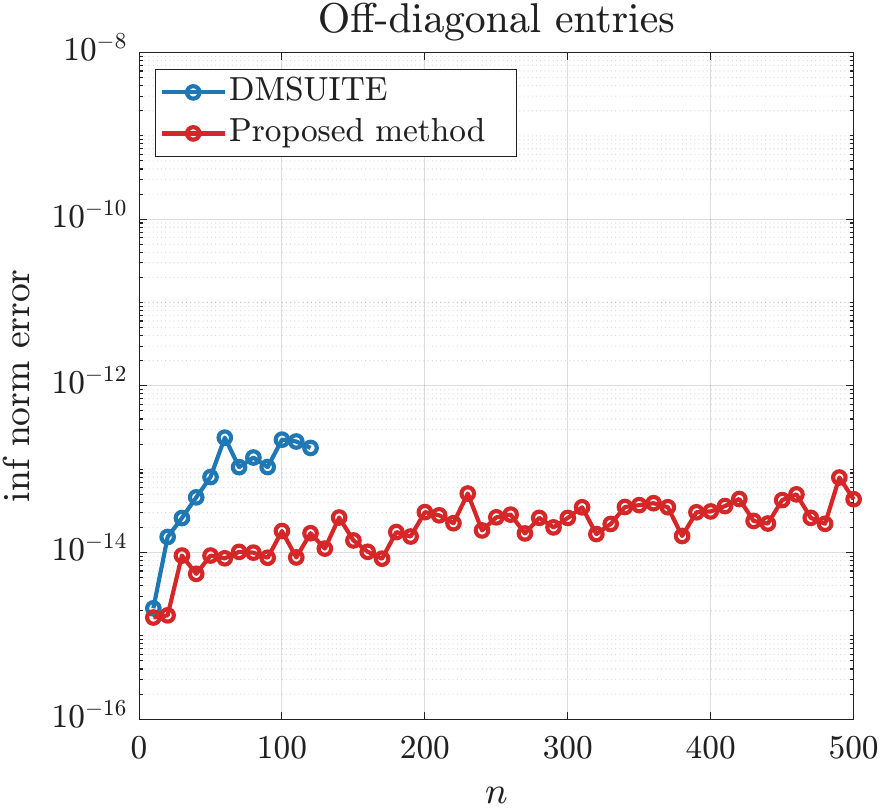}
    \caption{Relative error in the computation of a first-order Laguerre pseudospectral differentiation matrix of size $n \times n$, measured against a reference solution computed in high-precision arithmetic. 
    Left: Relative error in the diagonal entries. The error in the implementation by \texttt{DMSUITE}~\protect\citep{dmsuite} grows rapidly with $n$ compared to the proposed formula used in this work. Right: Relative error in the off-diagonal entries. The \texttt{DMSUITE} implementation exhibits numerical breakdown at approximately $n = 125$, whereas the proposed algorithm remains stable and maintains accuracy for significantly larger values of $n$.}
    \label{fig:breakdown}
\end{figure}

This limitation is reflected in commonly used implementations. Figure~\ref{fig:breakdown} illustrates the error incurred in the computation of a first-order Laguerre pseudospectral differentiation matrix as the number of collocation points increases. When using the widely adopted MATLAB software package \texttt{DMSUITE}~\citep{dmsuite}, numerical breakdown is observed in the off-diagonal entries of this matrix at approximately $n = 125$, accompanied by rapidly growing error in the diagonal entries. In contrast, the approach proposed in this work reformulates the expressions for the matrix entries and computes all quantities using a numerically stable algorithm, resulting in improved accuracy and a substantial increase in robustness.

\begin{remark}\em The \texttt{DMSUITE} package was published in 2003. \citet{weideman} notes that ``Almost all the codes will break down if the order of the matrix becomes too large; in the case of \texttt{lagdif.m} it is around $100 \times 100$ or so and larger for the others.'' and remarks that ``If you need to use matrices this large, it begs the question whether you should be using this particular pseudospectral method at all.'' This viewpoint reflects the computational context of the time; today such matrix sizes are no longer considered significantly large and numerical difficulties at this scale need not preclude the use of such methods.
\end{remark}

The purpose of this paper is to present a comprehensive approach for constructing Laguerre pseudospectral differentiation matrices that are stable, accurate, and practical for large-degree discretisations. The core contribution is a reformulation of the matrix entries that avoids catastrophic overflow, underflow, and loss of significance inherent in standard constructions. We demonstrate an improved formula for the off-diagonal entries, accompanied by a stable procedure for computing them, and motivate the use of simple closed-form expressions for the diagonal entries. This yields a procedure capable of reliably generating accurate differentiation matrices for substantially larger numbers of collocation points than is possible with classical implementations. Numerical experiments, including comparisons with the widely used \texttt{DMSUITE} package, demonstrate the improved robustness and accuracy.

The remainder of this paper is organised as follows. In Section~\ref{sec:preliminaries}, we review some necessary background on Laguerre polynomials and pseudospectral differentiation matrices. Section~\ref{sec:instability} examines the sources of numerical instability in classical constructions of Laguerre differentiation matrices. The stable construction proposed in this work is presented in Section~\ref{sec:construction}, followed by numerical experiments in Section~\ref{sec:experiments}. Concluding remarks are given in Section~\ref{sec:conclusion}.

\section{Preliminaries}
\label{sec:preliminaries}

The primary focus of this work is the formulation of differentiation matrices for Laguerre pseudospectral methods. In such methods, the solution of a differential equation is approximated by enforcing the governing equations to hold exactly at a finite set of prescribed nodes, known as \emph{collocation points}. Differentiation matrices are a key ingredient in these methods as they act on the function values at the collocation points and return derivative values at those points. To construct such matrices, we first discuss the underlying Laguerre polynomials, establish a choice of collocation nodes, and define the interpolation scheme used to build the matrix.

\subsection{Laguerre polynomials and functions}

\begin{figure}[H]
    \centering
    \includegraphics[height=5.5cm]{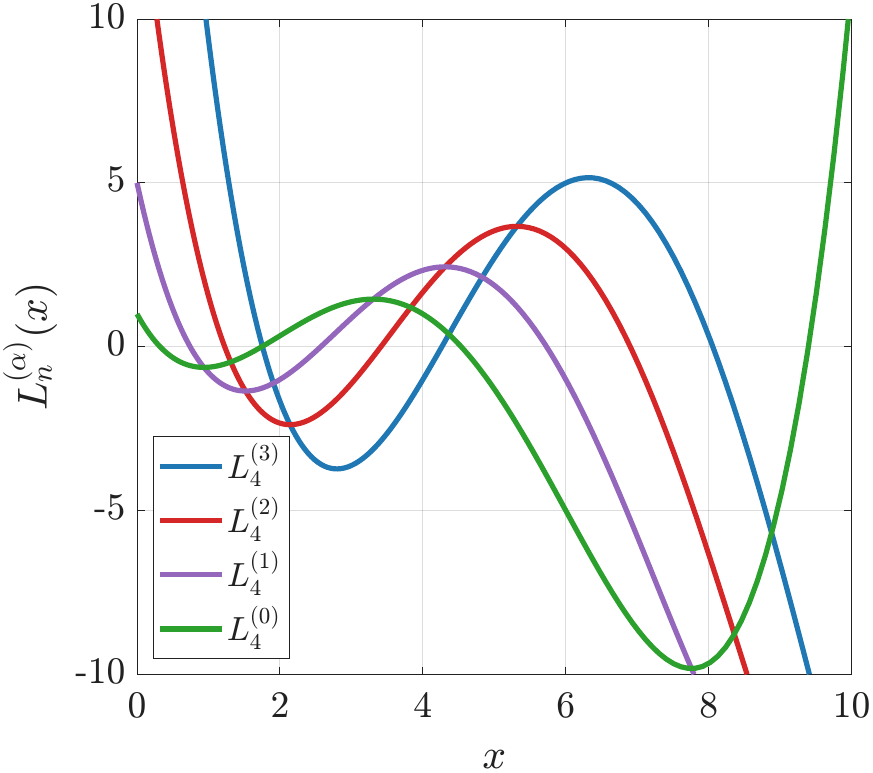}
    \hspace{15pt}
    \includegraphics[height=5.5cm]{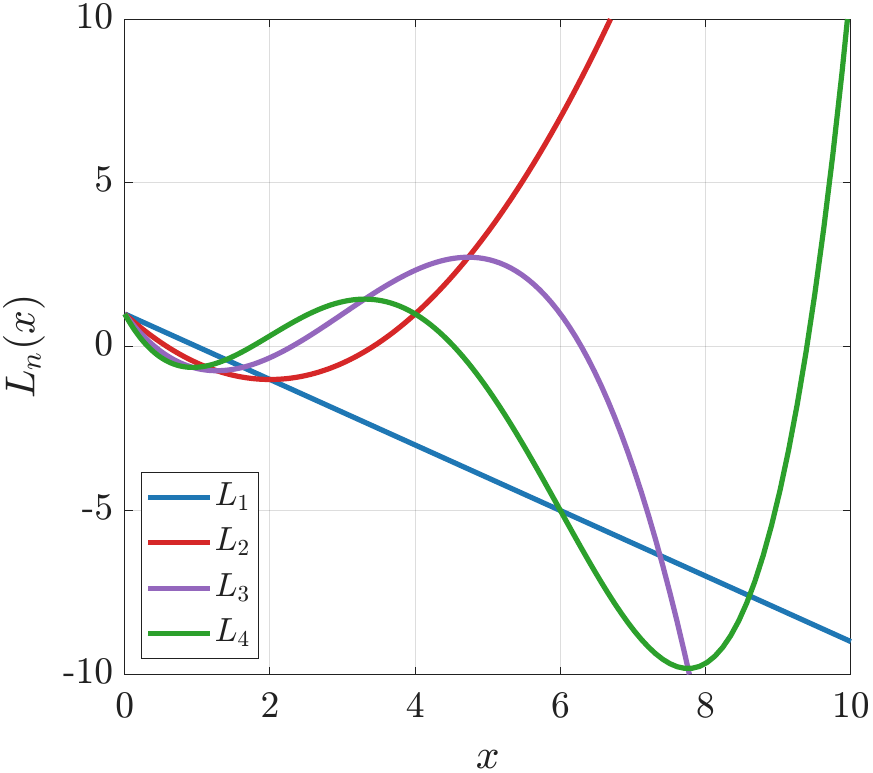}
    \caption{Left: The generalised Laguerre polynomials for $\alpha = 0, 1, 2, 3$ and $n = 4$. Right: The classical Laguerre polynomials ($\alpha = 0$) for $n = 1, 2, 3, 4$. For low degrees $n$ evaluated at small values $x$, the Laguerre polynomials are well-behaved. This is in contrast to the huge oscillations that occur for high degrees at large values of $x$, which is demonstrated in Figure~\ref{fig:lag_big}.}
    \label{fig:lagpoly}
\end{figure}

The \emph{generalised Laguerre polynomials} $L_n\supa(x)$ (also called the \emph{associated Laguerre polynomials}), named after Edmond Laguerre (1834-1886), are defined as the non-trivial solutions to the differential equation~\citep{szego1967}
\begin{equation}
    x y^{\prime \prime}(x) + (\alpha+1 - x) y^{\prime}(x) + n y(x) = 0 \,, 
\end{equation}
where $\alpha > -1$ and $n$ is a non-negative integer. These polynomials are orthogonal on $[0,\infty)$ with respect to the weight $x^\alpha e^{-x}$~\citep[\S 18.3]{nist2010}
\begin{equation}
    \int_0^\infty x^\alpha e^{-x} \, L_m\supa(x)\,L_n\supa(x)\,dx 
    = \frac{\Gamma(n+\alpha+1)}{n!}\,\delta_{m,n}\,, 
\end{equation}
where $\Gamma$ is the gamma function and $\delta_{m,n}$ is the Kronecker delta. An explicit representation is given by the Rodrigues formula~\citep[{(18.5.5)}]{nist2010}
\begin{equation}
    L_n\supa(x) 
    = \frac{1}{n!}x^{-\alpha} e^{x} \frac{d^n}{dx^n}\!\left(e^{-x} x^{n+\alpha}\right).
\end{equation}
The three-term recurrence relation derived from the differential equation is given by~\citep[{(18.9.1)}]{nist2010}
\begin{equation}
    (n+1) L_{n+1}\supa(x) 
    = (2n + \alpha + 1 - x) L_n\supa(x) - (n+\alpha)L_{n-1}\supa(x)\,, 
    \label{eq:rr}
\end{equation}
with initial terms $L_{0}\supa(x) = 1$ and $L_{1}\supa(x) = 1 + \alpha - x$. Numerical evaluation via this recurrence is the standard procedure and carries a cost proportional to $n$ per evaluation point. The first few generalized polynomials are visualised on the left in Figure~{\ref{fig:lagpoly}} for $\alpha = 0,1,2,3$ and $n=4$. In the case $\alpha = 0$, one recovers the \emph{classical Laguerre polynomials}, denoted $L_n(x)$. These are depicted on the right in Figure~{\ref{fig:lagpoly}} for the first few degrees on the interval $[0,10]$. 

\begin{figure}[H]
    \centering
    \includegraphics[width=6.4cm]{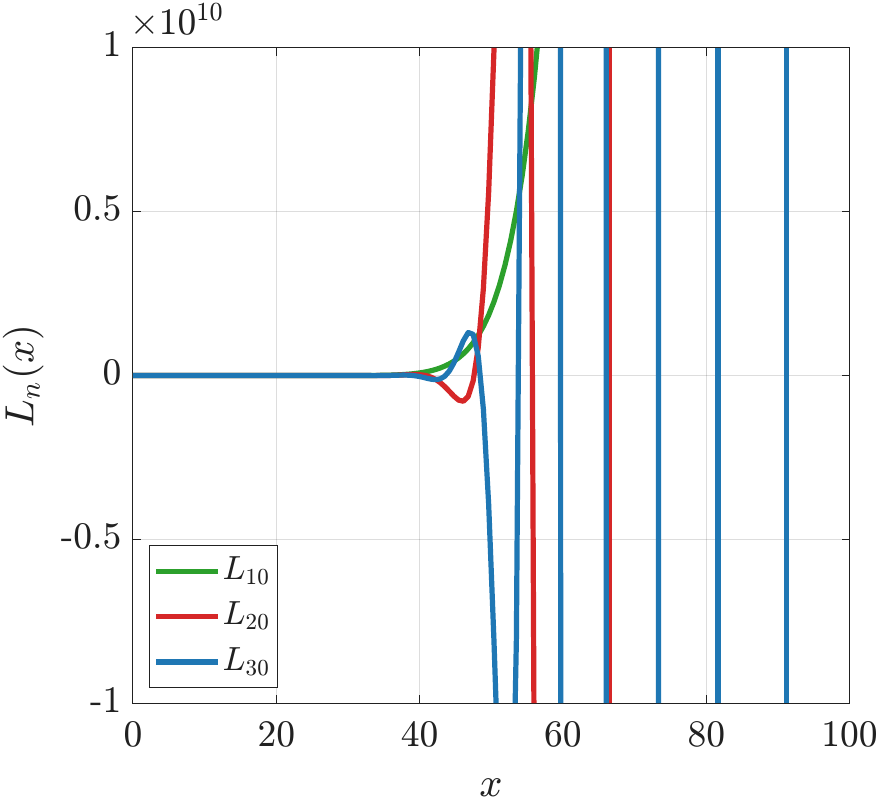}
    \hspace{10pt}
    \includegraphics[width=6.4cm]{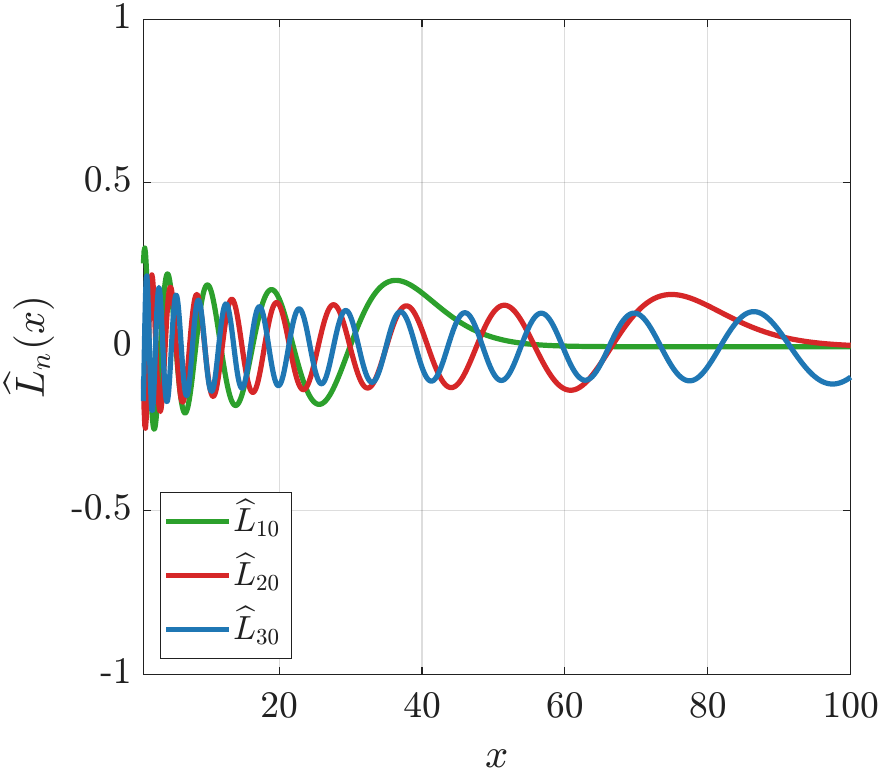}
    \caption{The classical Laguerre polynomials (left) and corresponding Laguerre functions (right) for degrees $n = 10, 20, 30$ as $x$ increases. The Laguerre polynomials exhibit severe oscillations and rapid growth at moderate degrees, whereas the inclusion of the exponential factor in the Laguerre functions results in well-scaled behaviour and significantly smaller magnitudes at large $x$.}
    \label{fig:lag_big}
\end{figure}

For small arguments $x$, the evaluation of Laguerre polynomials is generally unproblematic. However, even at moderate degrees, the polynomial grows rapidly as $x$ increases. As shown on the left of Figure~\ref{fig:lag_big}, the classical polynomials with degrees $n=10,20,30$ exhibit severe oscillations on $[0,100]$, with magnitudes exceeding $10^{10}$. Such growth can pose severe computational difficulties when using standard floating-point arithmetic. A common resolution is to work not directly with the polynomials, but instead with the associated \emph{Laguerre functions}, defined by
\begin{equation}
    \widehat{L}_n\supa(x) = e^{-x/2} L_n\supa(x).
    \label{eq:lagfunc}
\end{equation}
The exponential factor suppresses the growth of $L_n\supa(x)$ and ensures decay as $x \to \infty$. In fact, the Laguerre functions satisfy the bound $|\widehat{L}_n\supa(x)| \leq 1$ for all $x \geq 0$~\citep[(7.23)]{shen2011}. Furthermore, they maintain orthogonality, now with respect to the weight $x^{\scriptscriptstyle \alpha}$, share the same roots as the polynomials, and satisfy the corresponding differential equation~\citep[{(7.24)}]{shen2011}: 
\begin{equation} 
    x \,y^{\prime \prime}(x) + (\alpha+1)\,y^{\prime}(x) 
    + \left(n + \tfrac{1}{2}(\alpha+1) -\tfrac{1}{4}x \right)\,y(x) = 0.
\end{equation}
These functions may be generated using the same three-term recurrence relation as in~\eqref{eq:rr}, but with the initial terms scaled by the exponential weight, namely $\widehat{L}_0(x) = e^{-x/2}$ and $\widehat{L}_1(x) = (1 + \alpha - x)\,e^{-x/2}$. The resulting well-scaled behaviour is illustrated on the right in Figure~\ref{fig:lag_big}.

\subsection{Collocation points}
\label{subsec:lag_pseudo}

The roots of the Laguerre polynomials $L_n(x)$ provide a natural choice of collocation points for Laguerre pseudospectral methods. These points, known as the Laguerre--Gauss nodes, are also those used in Laguerre--Gauss quadrature~\citep[\S 7.1.2]{shen2011}. In many applications it is advantageous to include the endpoint $x=0$ explicitly in the collocation grid~\citep{weideman2000}. This can be achieved in two common ways. The first is to augment the Laguerre--Gauss nodes with the additional point $x_0=0$. The second is to use the Laguerre--Gauss--Radau nodes, which include the endpoint by construction.\footnote{This is analogous to the Chebyshev--Gauss--Lobatto nodes used in Chebyshev collocation, where both endpoints of the finite interval $[-1,1]$ are included.} These Radau nodes consist of $x_0=0$ together with the zeros of $L_n\supo(x)$~\citep[\S 3.2]{shen2009}.

Each set of nodes $\{x_j\}$ produces a distinct \emph{nodal polynomial} of the form
\begin{equation}
    \pi(x) = \prod_{\substack{j}} ( {x - x_j}).
    \label{eq:nodal_poly_prod}
\end{equation}
For the three distributions considered here, the nodal polynomial may alternatively be expressed in the form
\begin{equation}
\pi(x) = C_n \, a(x)\, L_n\supa(x)\,, 
\label{eq:nodal_poly}
\end{equation}
where $C_n = (-1)^n/n!$ denotes the leading coefficient of the degree-$n$ Laguerre polynomial, and the parameter $\alpha$ is given in Table~\ref{tab:laguerre_nodes}. The coefficient $a(x)$ is either $1$ or $x$, depending on whether the distribution includes an augmented node at the origin. Each nodal polynomial produces a corresponding distinct differentiation matrix.

\begin{table}[ht]
\centering
\renewcommand{\arraystretch}{1.4}
\caption{Laguerre node distributions of size $n+1$ and their associated parameters.}
\label{tab:laguerre_nodes}
\begin{tabular}{lll}
\hline
\\[-3ex]
Node distribution & Definition of nodes & Nodal polynomial parameters \\[6pt]
\hline
\\[-3ex]
Laguerre--Gauss 
& Zeros of $L_{n+1}(x)$ 
& $\alpha = 0$, \quad $a(x)= 1$ \\[6pt]

Augmented Laguerre--Gauss 
& $x=0$ and zeros of $L_n(x)$ 
& $\alpha = 0$, \quad $a(x)=x $ \\[6pt]

Laguerre--Gauss--Radau 
& $x=0$ and zeros of $L_n\supo(x)$ 
& $\alpha = 1$, \quad $a(x)=x$ \\[6pt]
\hline
\\[-4ex]
\end{tabular}
\end{table}

Because the standard Laguerre--Gauss nodes do not include the endpoint $x=0$, enforcing boundary conditions in the resulting collocation scheme requires some additional care. When the boundary is excluded, conditions at $x=0$ can still be imposed in two straightforward ways. One approach is to evaluate the global interpolating polynomial at $x=0$ and enforce this as an algebraic constraint, typically by replacing one of the interior collocation equations in the resulting linear system~\citep[\S 6.4]{boyd2001}. Alternatively, if the boundary conditions are linear, one may then employ ``basis recombination'', in which appropriate linear combinations of the original basis functions are constructed so that the resulting basis functions individually satisfy the boundary conditions~\citep[\S 6.5]{boyd2001}.

\subsection{Computation of zeros}

\begin{figure}[H]
    \centering
    \includegraphics[height = 6cm]{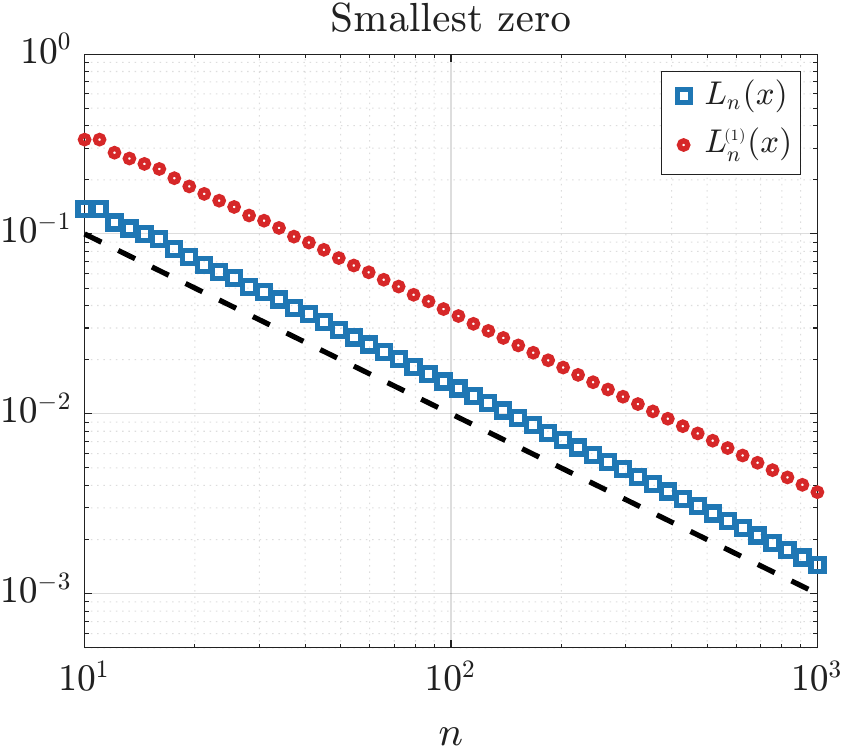}
    \hfill
    \includegraphics[height = 6cm]{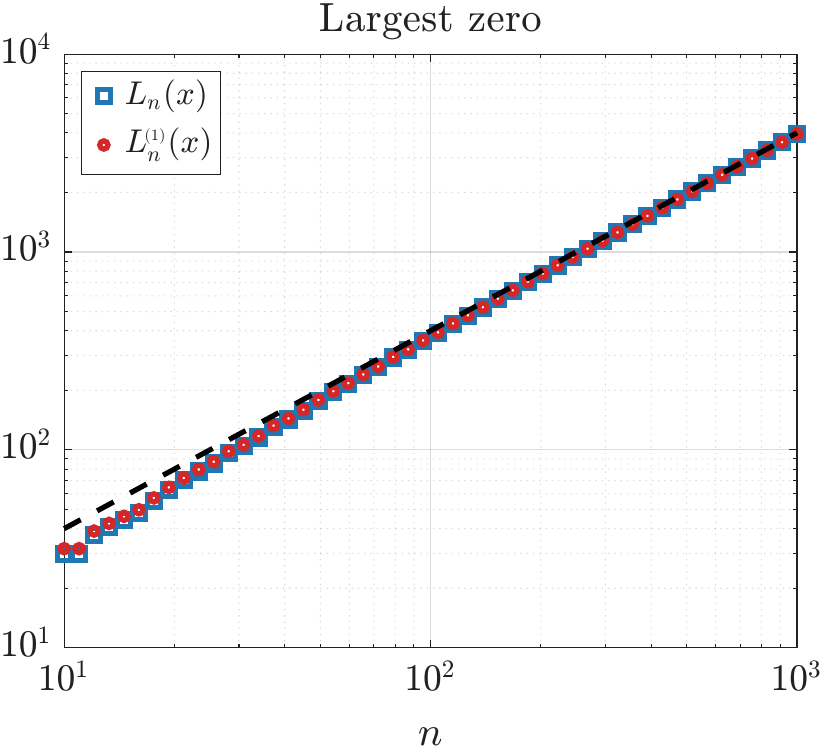}
    \caption{The smallest (left) and largest (right) zeros of $L_n(x)$ and $L_n\supo(x)$ as $n$ increases. The largest roots grow linearly like $4n$, while the smallest roots decay like $1/n$, as indicated by the black dashed reference lines.}
    \label{fig:roots}
\end{figure}

No explicit expressions are available for the zeros of Laguerre polynomials; however, asymptotic bounds provide useful estimates for large \(n\). In particular, the smallest and largest zeros of \(L_n\supa(x)\) satisfy the inequalities~\citep[{(6.31.7)} \& {(6.31.12)}]{szego1967}
\begin{equation}
    x_{1} > \frac{c}{2n + \alpha + 1 }, \qquad x_{n} < 2n + \alpha + 1 + \sqrt{(2n + \alpha + 1)^2 + \tfrac{1}{4} - \alpha^2}\,, 
    \label{eq:roots_asymp}
\end{equation}
where \(c > 0\) is a constant. The largest zero grows linearly like \(4n\), while the smallest zero decays on the order of \(1/n\)~\citep[\S 7.1.3]{shen2011}. These trends are illustrated in Figure~\ref{fig:roots}. Sharper bounds and refined asymptotic descriptions can be found in the work by~\citet{gatteschi2002}.

While explicit formulas are unavailable, efficient numerical methods for the accurate computation of the zeros have received significant attention. Early approaches relied on Newton iteration with sufficiently accurate initial guesses, however these methods tend to become unstable as $n$ increases. The Golub--Welsch algorithm~\citep{golub1969} reduced the problem to a symmetric tridiagonal eigenvalue problem that can be solved using fast eigensolvers (e.g.~\citep{gu1995}) in $O(n\log n)$ operations (or $O(n^2)$ operations when eigenvectors/quadrature weights are also computed). Over the past two decades, several algorithms achieving $O(n)$ complexity have been developed. One of the first was the iterative algorithm by~\citet{glaser2007}, later followed by the fourth-order fixed-point iteration method proposed by~\citet{gil2019}. In parallel, non-iterative, direct approaches based purely on asymptotic expansions have advanced considerably, with work by~\citet{opsomer2023}, and independently by~\citet{gil2018}, achieving double-precision accuracy even for moderate $n$ (e.g., $n > 100$). More recently, hybrid strategies combining both asymptotic-free iterative techniques and iteration-free asymptotics have also gained prominence~\citep{gil2025}. 

\subsection{Interpolation}

A fundamental ingredient in spectral collocation methods is interpolation. The exponential weighting typically associated with Laguerre methods is incorporated through \emph{weighted} polynomial interpolation~\citep{welfert1997,weideman2000}. The exponentially weighted interpolant of a function $f(x)$ at arbitrary nodes $\{x_j\}_{j=0}^n$ is given by
\begin{equation}
    f(x) \approx e^{-x/2}p_{n}(x) = \sum_{j=0}^n \frac{e^{-x/2}}{e^{-x_j/2}} \phi_j(x) f_j\,, 
    \label{eq:interpolant}
\end{equation}
where $f_j = f(x_j)$ and $\{\phi_j(x)\}$ are the Lagrange interpolating polynomials of degree $n$ defined by
\begin{equation}
\label{eq:interppoly}
    \phi_j(x) = \prod^n_{\substack{m=0 \\ m \neq j}} \left( \frac{x - x_m}{x_j - x_m}\right), \qquad j = 0,\ldots,n\,, 
\end{equation}
such that $\phi_j(x_k) = \delta_{j,k}$. The numerator of $\phi_j(x)$ is the {nodal polynomial} divided by $(x-x_j)$. The denominator is the \emph{barycentric weights} typically associated with the interpolation nodes, i.e.,
\begin{equation}
    v_j = \prod_{\substack{m = 0 \\ m \neq j}}^n \frac{1}{x_j - x_m}.
    \label{eq:vj_prod}
\end{equation}
Equivalently, in terms of the derivative of the nodal polynomial, the barycentric weights may also be expressed as~\citep[{p.~243}]{henrici1964}
\begin{equation}
    v_j = \frac{1}{\pi'(x_j)}.
    \label{eq:vj_der0}
\end{equation}
While the product form~\eqref{eq:vj_prod} provides a direct definition of the barycentric weights, it is numerically unstable at high polynomial degrees due to the accumulation of cancellation errors in the differences $(x_j - x_m)$~\citep{baltensperger2003}. In practice, the representation~\eqref{eq:vj_der0} is preferred, particularly when the nodes are related to classical orthogonal polynomials as in Eq.~\eqref{eq:nodal_poly}. The Lagrange interpolating polynomials in~\eqref{eq:interppoly} can be written in terms of the barycentric weights and nodal polynomial as:
\begin{equation}
    \phi_j(x) = \frac{\pi(x)}{x - x_j} v_j.
    \label{eq:phi}
\end{equation}

\subsection{Laguerre pseudospectral differentiation matrices}
\label{subsec:difmat}

In collocation methods, differentiation is performed via pseudospectral differentiation matrices. These matrices are generated by differentiating the interpolant~\eqref{eq:interpolant} and evaluating at the collocation nodes $\{x_k\}$:
\[
f^{(\ell)}(x_k) := \sum_{j=0}^{n} \frac{d^{\ell}}{dx^{\ell}} \left[  \frac{e^{-x/2}}{e^{-x_j/2}} \phi_j(x)  \right]_{x = x_k} f(x_j)\,, 
\qquad k = 0, \dots, n.
\]
The differentiation matrix of order $\ell$ may be represented by a matrix $D^{(\ell)}$, with entries~\citep{weideman2000}:
\begin{equation}
D^{(\ell)}_{k,j} := \frac{d^{\ell}}{dx^{\ell}} \left[\frac{e^{-x/2}}{e^{-x_j/2}} \phi_j(x)\right]_{x = x_k}.
\label{eq:derdifmat}
\end{equation}
Differentiation is thus reduced to the matrix-vector product $\bm{f}^{(\ell)} \approx D^{(\ell)} \bm{f}$, where $\bm{f} = [f(x_0),\, f(x_1), \, \ldots, \,f(x_n)]^\top $ is the vector of function values at the collocation points and $\bm{f}^{(\ell)}$ contains the approximate $\ell^{\text{th}}$ derivatives at those nodes. For $\ell = 1$, the entries in the first-order differentiation matrix $D$ can be expanded via the product rule:
\begin{equation}
D_{k,j} = -\frac{1}{2}\frac{e^{-x_k/2}}{e^{-x_j/2}} \phi_j(x_k) + \frac{e^{-x_k/2}}{e^{-x_j/2}} \phi^{\prime}_j(x_k).
\end{equation}
Separating \emph{diagonal} ($j = k$) and \emph{off-diagonal} ($j \neq k$) entries yields the formulas: 
\begin{equation}
D_{k,k}  = -\frac{1}{2} + \phi^{\prime}_k(x_k),  \qquad D_{k,j} =\frac{e^{-x_k/2}}{e^{-x_j/2}} \phi^{\prime}_j(x_k)\,, 
\label{eq:difmat_entries}
\end{equation}
for $j,k = 0, \ldots, n,\; j \neq k$. 

Differentiating the interpolating polynomials in~\eqref{eq:phi} and evaluating at $x=x_k$ for $j \neq k$ produces
\begin{equation}
    \phi^{\prime}_j(x_k) = \frac{v_j/v_k}{x_k-x_j}\,, 
    \label{eq:phider}
\end{equation}
where the identity in~\eqref{eq:vj_der0} has been used to simplify the expression. Substituting into~\eqref{eq:difmat_entries} gives the off-diagonal entries of $D$ as:
\begin{equation}
    D_{k,j} =  \frac{c_k/c_j}{x_k-x_j} \quad \text{where}\quad c_j = \frac{e^{-x_j/2}}{v_j}, \quad j \neq k.
    \label{eq:Dkj}
\end{equation}
It is standard practice in collocation methods to compute the diagonal entries using the off-diagonal entries---a step commonly referred to as the {`negative sum trick'}~\citep{baltensperger2003}.  In weighted polynomial interpolation, the negative sum trick enforces that the differentiation matrix exactly differentiates the weight function at the collocation points~\citep{weideman2000}, e.g.,
\begin{equation}
    D e^{-\bm{x}/2} = -\frac{1}{2} e^{-\bm{x}/2} \quad \implies \quad \sum_{j = 0}^n D_{k,j} \, e^{-x_j/2} = -\frac{1}{2} e^{-x_k/2} \quad \text{for each $k$}.
\end{equation}
This condition leads to the diagonal entries expressed using the sum of the off-diagonal entries:
\begin{equation}
    D_{k,k}  =-\frac{1}{2} - \sum_{\substack{j = 0\\ j \neq k}}^n \frac{e^{-x_j/2}}{e^{-x_k/2}} \, D_{k,j}.
    \label{eq:Djj}
\end{equation}
While in unweighted polynomial interpolation, the negative sum trick provides a simple and accurate way to compute the diagonals, in the weighted setting it becomes unstable for large $n$. This is due to the exponential factors in~\eqref{eq:Djj} that amplify round-off errors in the off-diagonals, particularly for the first columns of the last few rows. Instead, an alternative formula may be obtained by substituting an expression for $\phi_k'(x_k)$ into~\eqref{eq:difmat_entries}. To this end, we differentiate the expression in~\eqref{eq:phi} for $j = k$ twice to obtain
\begin{equation}
    (x-x_k)\phi_k^{\prime \prime}(x)+2\phi_k^{\prime}(x)=v_k\,\pi^{\prime \prime}(x).
\end{equation}
Evaluating at $x=x_k$ and utilising the identity in~\eqref{eq:vj_der0} yields
\begin{equation}
    \phi_k^{\prime}(x_k)
    = \frac{1}{2}v_k\,\pi^{\prime \prime}(x_k)
    = \frac{\pi^{\prime \prime}(x_k)}{2\pi^{\prime}(x_k)}.
\end{equation}
A direct differentiation of the nodal polynomial in~\eqref{eq:nodal_poly_prod} shows that
\begin{align*}
    \pi^{\prime \prime}(x_k) &= 2\sum_{\substack{i = 0 \\ i \neq k}}^n  \prod_{\substack{j = 0 \\ j \neq i,k}}^n (x_k-x_j)
    = 2\sum_{\substack{i = 0 \\ i \neq k}}^n   \frac{1}{x_k-x_i}\prod_{\substack{j = 0 \\ j \neq k}}^n (x_k-x_j) 
    = 2\,\pi^{\prime}(x_k)\sum_{\substack{i = 0 \\ i \neq k}}^n   \frac{1}{x_k-x_i}\,, 
\end{align*}
and hence
\begin{equation}
    \phi_k^{\prime}(x_k)
    = \sum_{\substack{i=0 \\ i\neq k}}^n \frac{1}{x_k-x_i}.
\end{equation}
This returns
\begin{equation}
    D_{k,k}  = -\frac{1}{2} - \sum_{\substack{i=0 \\ i\neq k}}^n \frac{1}{x_k-x_i}\,, 
    \label{eq:Djj_simplified}
\end{equation}
which is the expression more commonly used in practice.

For higher order differentiation matrices, \citet{welfert1997} derived the following recursion for the off-diagonal entries
\begin{equation}
D^{(\ell)}_{k,j} =
\begin{cases}
\dfrac{\ell}{x_k - x_j}
\left(
\dfrac{c_k}{c_j} D^{(\ell-1)}_{k,k}
- D^{(\ell-1)}_{k,j}
\right)\,,  & j \neq k\,,  \\[2ex]
\left(-\dfrac{1}{2}\right)^\ell
-\displaystyle\sum_{\substack{j = 0 \\ j \neq k}}^{n}
\dfrac{e^{-x_j/2}}{e^{-x_k/2}} \, D^{(\ell)}_{k,j}\,, 
& j = k\,, 
\end{cases}
\label{eq:welfertrecursion}
\end{equation}
for $\ell = 1,2,\dots,$ with $D^{(0)}$ the identity matrix. The diagonal entries result directly from the negative sum trick applied to higher-order derivatives of the weight function~\citep{baltensperger2003}. While not typically used in practice, it provides a compact formulation. \citet{welfert1997} also proposed a recursion for the diagonal entries, however it requires the introduction of reduced differentiation matrices $\widehat{D}$, formed by removing a single node from the grid and constructing differentiation matrices on the remaining nodes. In Section~\ref{sec:construction}, we present an alternative formula specifically for second-order differentiation matrices that avoids this additional complexity.

Higher-order differentiation matrices can, in principle, also be obtained by taking successive powers of the first-derivative matrix. However, computing powers of a dense matrix requires $\mathcal{O}(n^3)$ operations, while the recursive algorithm of Welfert requires only $\mathcal{O}(n^2)$ operations. The recursive method also reduces the accumulation of round-off errors that tend to occur during repeated matrix multiplications~\citep{weideman2000}.

\subsection{Software for Laguerre pseudospectral differentiation matrices}
\label{subsec:software}

The landscape of software for spectral and pseudospectral methods has expanded significantly over the past two decades. Robust open-source packages are now available in multiple programming languages for solving ordinary and partial differential equations using global polynomial and trigonometric bases. In Python, frameworks such as \texttt{Dedalus}~\citep{dedalus} and \texttt{Shenfun}~\citep{shenfun} provide powerful environments for solving partial differential equations using spectral Galerkin and Tau formulations. Similarly, the Julia language offers the package \texttt{ApproxFun}~\citep{approxfun}, which leverages adaptive approximation with infinite-dimensional operators to solve differential equations using spectral techniques. In the GNU Octave ecosystem, the \texttt{SPSMAT}~\citep{spsmat} package supplies a general toolbox for spectral and pseudospectral discretisation.

Beyond these general-purpose solver frameworks, a number of more specialised packages target particular polynomial families and collocation strategies. For Chebyshev and Legendre methods in particular, the available software is extensive and mature, with highly optimised algorithms for node generation, quadrature weights, and differentiation matrices implemented across several environments. For example, \texttt{Chebfun}~\citep{chebfun} in MATLAB provides a comprehensive system centred on Chebyshev technology for function approximation and the numerical solution of differential equations.

By contrast, Laguerre spectral methods have received comparatively limited dedicated software support. The most prominent and widely used implementation of Laguerre pseudospectral differentiation matrices remains the Laguerre routine included in the Differentiation Matrix Suite (DMSUITE) of~\citet{dmsuite}. Originally developed for MATLAB and GNU Octave, with ports to Python~\citep{python_dmsuite} and Julia~\citep{julia_dmsuite} following later, \texttt{DMSUITE} provides a general-purpose framework for constructing pseudospectral differentiation matrices for several polynomial families, including Chebyshev, Legendre, Hermite, Fourier, sinc, and Laguerre. However, the Laguerre functionality implemented via the \texttt{lagdif} routine is not specialised for Laguerre methods and relies on classical algorithms that are not suitable for high-degree computations. In particular, the method is not robust for large numbers of collocation points and is sensitive to round-off error, as illustrated in Figure~\ref{fig:breakdown}. 

Another notable MATLAB implementation is the Laguerre spectral/pseudospectral library of~\citet{jewell2026}, which computes differentiation matrices based on scaled Laguerre functions rather than the standard (unscaled) Laguerre basis. While effective for certain applications, this formulation produces operators that differ from those associated with the classical Laguerre pseudospectral collocation method and is therefore not directly compatible with standard Laguerre discretisations. Other MATLAB packages exist that generate pseudospectral differentiation matrices on arbitrary grids in a manner similar to \texttt{DMSUITE}, e.g.~\citep{vonwinckel2026}. However, these general-purpose tools are not formulated within the weighted interpolation framework required for Laguerre methods. As a result, their direct application to Laguerre collocation leads to severe numerical overflow for even moderate polynomial degrees, rendering them unsuitable for reliable computation on semi-infinite domains.

\section{Instability of classical Laguerre differentiation matrices}
\label{sec:instability}

Despite their theoretical appeal, classical Laguerre collocation methods based on differentiation matrices are known to suffer from severe numerical instabilities as the number of nodes increases~\citep{funaro1990,shen2000,huang2024}. We now discuss the causes of these instabilities. 

\subsection{Sources of numerical instability}

\begin{figure}[ht]
    \centering
    \includegraphics[height=6cm]{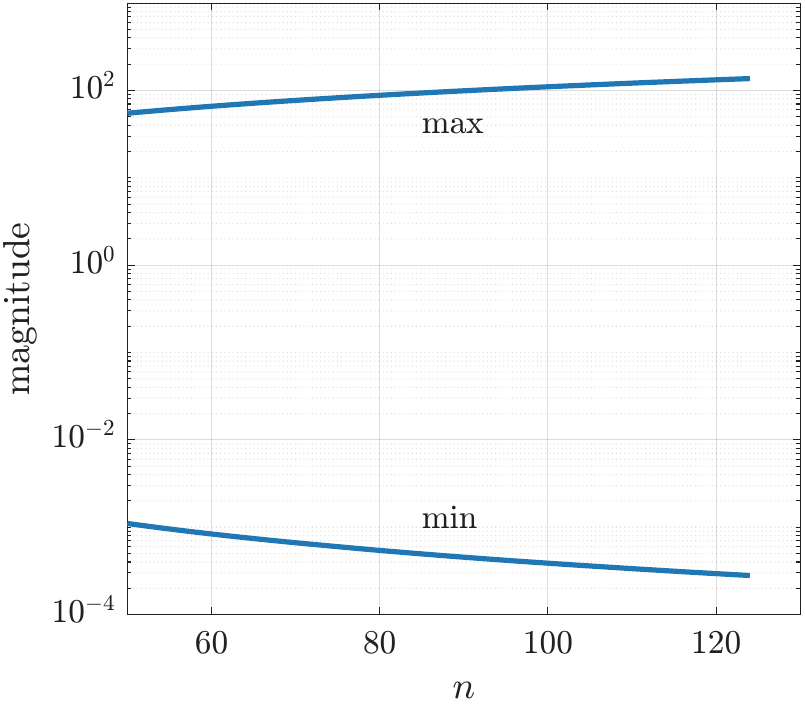}
    \hfill
    \includegraphics[height=6cm]{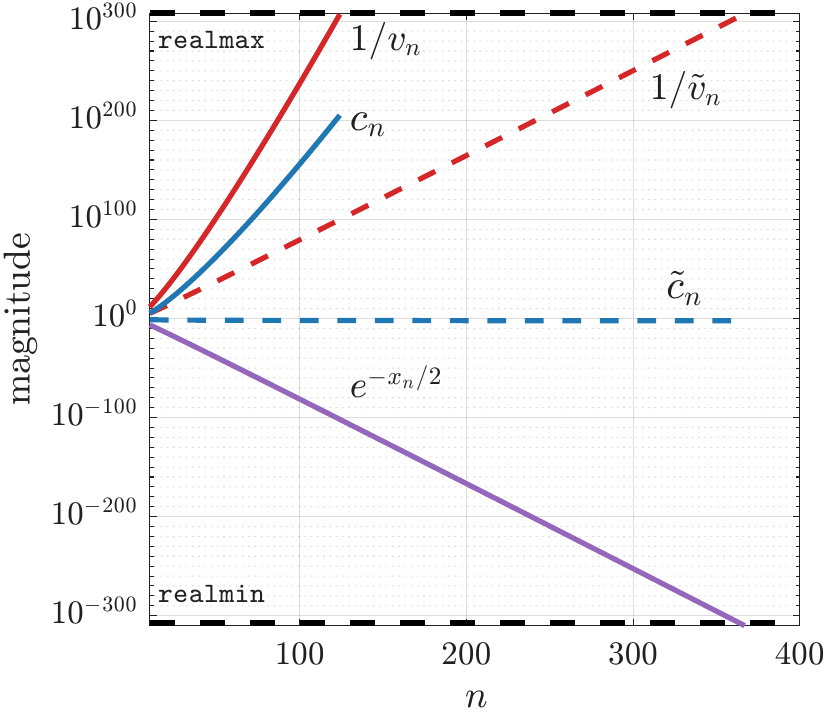}
    \caption{Left: The minimum and maximum magnitude entries in the first-order differentiation matrix constructed using \texttt{DMSUITE}~\protect\citep{dmsuite}. While these values are not significantly large or small, the computation fails at due to the intermediate quantities needed to compute the entries exceeding the \texttt{realmax} and \texttt{realmin} thresholds in IEEE double precision arithmetic. Right: The growth and decay of the factors defining $c_n$ (unscaled) and $\tilde c_n$ (scaled) at the largest node $x_n$. It shows $1/v_n$ grows rapidly and exceeds double-precision limits, causing the failure seen on the left. In contrast, the scaled quantity $1/\tilde{v}_n$ grows more moderately and balances the decay of $e^{-x_n/2}$, resulting in a well-scaled product $\tilde{c}_n$ up to around $n \approx 383$, when the components still individually overflow and underflow.}
    \label{fig:minmax}
\end{figure}

Numerical instability in Laguerre pseudospectral methods arises from three closely related mechanisms: (i) the rapid growth of Laguerre polynomials and their derivatives, (ii) exponential decay in the associated Laguerre functions, and (iii) the cancellation of extreme intermediate quantities in standard differentiation matrix formulas. Importantly, the final discrete operators are typically well-scaled. As illustrated in Figure~\ref{fig:minmax}, the maximum and minimum magnitude entries in the first-order Laguerre differentiation matrix produced by~\texttt{DMSUITE} are not especially large or small, yet numerical breakdown occurs at $n \approx 125$ due to overflow and underflow. The difficulty does not lie with the magnitude of the final entries, but rather with the conventional construction of these matrices, which requires evaluating intermediate terms that exceed the representation limit of double-precision arithmetic.

We first consider the behaviour of the Laguerre polynomials themselves. In pseudospectral discretisations, high-degree polynomials and their derivatives must be evaluated at increasingly large nodes. The magnitude of $L_n\supa(x)$ grows rapidly as $n$ and $x$ increase, and for $n$ in the several hundreds its values can exceed the overflow threshold of IEEE double precision. Direct polynomial-based formulations therefore become unreliable for sufficiently large discretisations.

To suppress this growth, many formulations instead employ the Laguerre functions $\widehat{L}_n\supa$ defined in~\eqref{eq:lagfunc}, which incorporate the weight $e^{-x/2}$. Although this removes polynomial blow-up in the final quantity, it introduces the complementary issue of exponential underflow. In double precision, $e^{-x/2}$ underflows once $x \gtrsim 1492$, and the largest root of $L_n\supa$ satisfies $x_{n} \approx 4n$, thus exponential underflow becomes unavoidable once $n \gtrsim 373$. 

These same mechanisms appear in the construction of Laguerre pseudospectral differentiation matrices, whose off-diagonal entries involve ratios $c_j/c_k$ with
\begin{equation} 
    c_j = \frac{e^{-x_j/2}}{v_j}\,, 
    \label{eq:cj_def}
\end{equation}
where $v_j$ denotes the barycentric weights. The exponential decay in the numerator is countered by rapid growth in $1/v_j$. Standard implementations, such as \texttt{DMSUITE}, evaluate the numerator and denominator separately before multiplying.  For large degrees this procedure fails because the individual quantities exceed the bounds of double-precision arithmetic, even though their product does not. This behaviour is illustrated in Figure~\ref{fig:minmax} for the largest node $x_n$, where $1/v_n$ exceeds the \texttt{realmax} threshold.

Since the differentiation matrix depends only on ratios of these quantities, $v_j$ may be rescaled by any factor depending solely on $n$. This freedom allows for a more stable representation. Using the representation of the barycentric weights in terms of the nodal polynomial in~\eqref{eq:vj_der0}, the weights corresponding to $n$ Laguerre--Gauss nodes (Table~\ref{tab:laguerre_nodes}) can be written as
\begin{equation}
    v_j = \frac{1}{\pi'(x_j)} 
        = \frac{1}{C_n L_n'(x_j)}\,, 
    \label{eq:vj_der}
\end{equation}
where $C_n = (-1)^n/n!$. Compared to the product formulation in~\eqref{eq:vj_prod}, this expression has the advantage of explicitly expressing the rapidly growing factor $1/C_n$. Since this factor cancels in ratios, it can be omitted. We therefore define the scaled quantities
\begin{equation}
    \tilde v_j = \frac{1}{L_n'(x_j)}\,,  \qquad 
    \tilde c_j = \frac{e^{-x_j/2}}{\tilde v_j}.
    \label{eq:cj_scaled}
\end{equation}
The ratios $\tilde c_j/\tilde c_k$ yield the same differentiation matrix as $c_j/c_k$, but the magnitude of intermediate quantities is reduced. As shown in Figure~\ref{fig:minmax}, the growth of $1/\tilde v_n$ is much slower than that of $1/v_n$, but more importantly, the growth is now at a rate similar to the decay of $e^{-x_n/2}$. As a result, the product $\tilde c_n$ remains well-scaled up to $n \approx 383$, where it is now the underflow in the exponential term that causes failure. 

This demonstrates that the above rescaling delays the onset of numerical breakdown, but does not eliminate it. To avoid failure entirely, it is necessary to compute $\tilde c_j$ directly, rather than as the product of two separately over- and underflowing terms.

\subsection{Existing methods for avoiding instabilities}

Several approaches have been proposed to address the numerical instabilities that arise due to the rapid growth of Laguerre polynomials. These methods primarily aim to stabilise the evaluation of Laguerre functions at high polynomial degrees. \citet{funaro1990} introduced \emph{scaled} Laguerre functions obtained by applying an $n$-dependent scaling to the classical Laguerre polynomials. This scaling significantly reduces the dynamic range of the basis functions and enables stable evaluation for larger values of $n$, however the scaled functions are not simple and increases implementation complexity, particularly for derivative computations. \citet{shen2000} proposed a sign and amplitude factorisation technique that evaluates Laguerre functions by separating their sign from a suitably scaled magnitude. More recently, \citet{huang2024} developed an adaptive stabilisation strategy in which the exponential weight is incorporated dynamically during the recurrence. Intermediate quantities are rescaled whenever their magnitude exceeds a prescribed threshold, and the remaining weight is applied at the end of the computation. While successful, the drawbacks are additional algorithmic complexity and suitable parameter selection.

While these works propose methods to mitigate sources of numerical instability in the evaluation of Laguerre polynomials and functions, they do not discuss or address the stable construction of Laguerre differentiation matrices in collocation settings. We aim to do this here.

\section{Stable construction of Laguerre differentiation matrices}
\label{sec:construction}

We propose a strategy for the stable computation of Laguerre pseudospectral differentiation matrices. Our approach treats the off-diagonal and diagonal entries separately. For the off-diagonal entries, we make two key observations. First, entries should be computed using scaled coefficients $\tilde{c}_j$ such as in~\eqref{eq:cj_scaled} rather than the unscaled $c_j$. Second, and perhaps most importantly, the quantities making up these coefficients should be computed simultaneously rather than in isolation, as separate evaluation leads to overflow and underflow. To facilitate this, we utilise an algorithm introduced by~\citet{glaser2007}, which enables the required values to be computed in a unified and stable manner.

For the diagonal entries, we employ a direct closed-form expression. Although such formulas exist, we provide a derivation in the Laguerre setting and offer justification for preferring this direct approach over alternatives used in literature. Together, these components yield a comprehensive and stable method for constructing Laguerre pseudospectral differentiation matrices at large $n$.

\subsection{Reformulated expression for off-diagonal entries}
\label{subsec:offdiag}

We now describe the scaled barycentric weights and corresponding coefficients for the distributions considered in terms of the derivatives of Laguerre polynomials. Differentiating the form of the nodal polynomial in~\eqref{eq:nodal_poly} via the product rule produces
\begin{equation}
    \pi^{\prime}(x) = C_n \, a^{\prime}(x)\,L_n\supa(x) + C_n\,a(x)\,L_n\supap(x)\,, 
\end{equation}
with $a(x)$ given in Table~\ref{tab:laguerre_nodes}. For the augmented and Radau distributions, $a(0) = 0$ and $a^{\prime}(0) = 1$, thus evaluating at the augmented node $x_0 = 0$ produces:
\begin{equation}
    \pi^{\prime}(0) = C_n \,L_n\supa(0) = C_n\binom{n+\alpha}{n}\,, 
\end{equation}
using the identity~\citep[{(18.6.1)}]{nist2010}. Alternatively, at any root $x_j$ of $L_n\supa(x)$, the first term is zero and the derivative simplifies to:
\begin{equation}
    \pi^{\prime}(x_j) = C_n \, a(x_j)\,L_n\supap(x_j).
\end{equation}
The identity in~\eqref{eq:vj_der0} gives the barycentric weights as the inverse of these derivatives. Hence the scaled barycentric weights corresponding to nonzero nodes are
\begin{equation}
    \tilde{v}_j = \frac{1}{a(x_j)\,L_n\supap(x_j)}\,,  \quad j = 1,\ldots,n,
\end{equation}
and $\tilde{v}_0 = 1/\binom{n+\alpha}{n}$ for the zero node. Scaling by $C_n$ is necessary to exclude the fast growing factorial from these expressions. The corresponding scaled coefficients $\tilde{c}_j$ are
\begin{equation}
    \tilde{c}_j
    = a(x_j)\, L_n\supap(x_j)e^{-x_j/2}\,, 
    \qquad j=1,2,\ldots,n,
\end{equation}
and $\tilde{c}_0 = \binom{n+\alpha}{n}$. The exponential factor can now be absorbed by the Laguerre polynomials by converting to the Laguerre functions $\widehat{L}_n\supa(x)$, since
\begin{equation}
\widehat{L}_n\supap(x_j)=\frac{d}{dx} \left[ e^{-x_j/2} L_n\supa(x_j) \right]_{x=x_j} = e^{-x_j/2} L_n\supap(x_j).
\end{equation}
This leaves:
\begin{equation}
    \tilde{c}_j
    = a(x_j)\, \widehat{L}_n\supap(x_j)\,, 
    \qquad j=1,2,\ldots,n.
    \label{eq:cj_new}
\end{equation}
To stably compute $\tilde{c}_j$ in this form, we propose a method for evaluating the derivatives of Laguerre functions at their roots that avoids numerical overflow or underflow. Our approach relies on a root-finding algorithm. In the following section, we introduce this algorithm and demonstrate how the required derivatives are naturally generated in a way that avoids underflow and overflow as part of the root-finding process.

\subsection{Glaser--Liu--Rokhlin algorithm}
\label{subsec:glr}

\citet{glaser2007} introduced a fast and stable algorithm for computing the roots of a function $y(x)$ satisfying a second-order ordinary differential equation of the form
\begin{equation}
p(x)\,y^{\prime \prime}(x) + q(x)\,y^{\prime}(x) + r(x)\,y(x) = 0\,, 
\end{equation}
where $p(x)$, $q(x)$, and $r(x)$ are polynomials of degree at most two. The roots are computed sequentially using a predictor-corrector approach. To advance from a known root $x_i$ to the next, the algorithm first obtains an initial approximation for $x_{i+1}$ by numerically solving a first-order auxiliary differential equation derived via the Pr\"{u}fer transform. Once this approximate root is found, Newton iteration is applied as a corrector step to determine the precise distance $h$ to the root, i.e., $x_{i+1} = x_i + h$. Each Newton update takes the form
\begin{equation}
h^{[j+1]} = h^{[j]} - \frac{y(x_i + h^{[j]})}{y^{\prime}(x_i + h^{[j]})}\,, 
\label{eq:newton}
\end{equation}
which requires evaluating the function and its derivative at the next root. Rather than evaluating $y$ and $y^{\prime}$ directly, the algorithm computes these quantities using local Taylor expansions about the known root $x_i$:
\begin{equation}
y(x_i+h) = \sum_{k=0}^m \frac{y^{(k)}(x_i)}{k!} h^k, 
\qquad
y^{\prime}(x_i+h) = \sum_{k=1}^m \frac{y^{(k)}(x_i)}{(k-1)!} h^{k-1}.
\label{eq:taylor}
\end{equation}
The higher-order derivatives $y^{(k)}(x_i)$ are obtained recursively from the differential equation satisfied by $y$. By repeatedly differentiating the governing equation and evaluating it at $x_i$, each derivative can be expressed in terms of lower-order derivatives. 

To initialise the procedure, an initial root $x_1$ is required. This is obtained by selecting a starting value $x_s < x_1$, from which the Pr\"{u}fer-based predictor and subsequent Newton refinement are applied to compute an accurate value of $x_1$. For orthogonal polynomials and related special functions, suitable choices of $x_s$ are available from known inequalities and asymptotic bounds on the smallest root~\citep{gautschi2004,szego1967} (for example, Eq.~\eqref{eq:roots_asymp} for the Laguerre polynomials). In addition to the location of the initial root, the values of $y$ and $y^{\prime}$ at $x_1$ must be available. These are typically computed using a scheme specific to the function $y$; for orthogonal polynomials it is the associated three-term recurrence relation.

This procedure, referred to as the Glaser--Liu--Rokhlin algorithm, offers two primary advantages in the present setting: 1) it produces accurate values for the roots, which serve as the collocation nodes in our pseudospectral method, and 2) it evaluates $\widehat{L}_n\supap(x)$ at these roots directly, without the explicit formation of exponentially large or small intermediate quantities. We apply the algorithm to the Laguerre functions $\widehat{L}_n\supa$ rather than the polynomials; while they share the same roots, the functions ensure that all derivative values remain well-scaled. The key benefit is that once the exponential weight is absorbed into the function definition at the initial root, all subsequent evaluations are performed via local Taylor expansions or the differential equation satisfied by $\widehat{L}_n\supa$. Consequently and importantly, the factor $e^{-x/2}$ is never explicitly evaluated at large values of $x$, thereby avoiding numerical overflow or underflow.

As noted by~\citet{glaser2007}, a specific complication arises when applying this algorithm to generate Laguerre roots. The differential equations satisfied by the Laguerre polynomials and functions exhibit a singularity at $x=0$, which causes the roots to cluster densely near the origin. Consequently, Taylor expansions centred at points close to zero may converge slowly or fail altogether. A practical remedy is to modify the algorithm for the first $k$ roots, where $k$ is chosen heuristically so that $x_k$ lies sufficiently far from the singularity.\footnote{In the implementation of the Glaser--Liu--Rokhlin algorithm in \texttt{lagpts.m} from the \texttt{Chebfun} package~\citep{chebfun} this heuristic value is taken as $k=20$.} For the roots $x_1,\ldots,x_k$, the evaluations of $y$ and $y^{\prime}$ required in~\eqref{eq:newton} are performed using the three-term recurrence relation satisfied by the Laguerre functions (or a modified variant for improved accuracy; see Section~\ref{subsec:modified}), rather than the Taylor expansions in~\eqref{eq:taylor}. This adaptation avoids convergence issues near the origin while still restricting direct evaluations to small values of $x$, where overflow and underflow are not a concern. Once $x_k$, together with $\widehat{L}_n\supa(x_k)$ and $\widehat{L}_n\supap(x_k)$, has been computed, all subsequent evaluations use the stable local expansions in~\eqref{eq:taylor}.

\subsection{Direct formula for diagonal entries}
\label{subsec:diagonals}

We now turn to the computation of the diagonal entries of the differentiation matrix. In Section~\ref{subsec:difmat}, the expression~\eqref{eq:Djj_simplified} was given as the form most commonly used in practice. In this section, we derive a simple and explicit formula for these entries.

Direct formulas for the diagonal entries exist for many polynomial families, including Laguerre (see, e.g.,~\citet[\S 7.1.5]{shen2011}), but they are often avoided in practice due to numerical instability arising from cancellation. For example, in the Chebyshev case, such formulas involve factors of the form $1/(1-x^2)$ and therefore suffer significant loss of precision near the endpoints $x=\pm 1$. As a result, the negative sum trick is generally regarded as the better approach. However, we shall see below that, in the Laguerre setting (at least for first-order differentiation matrices), cancellation error is not an issue. The direct formula remains numerically stable across the entire domain, and there is no inherent reason to avoid its use in computing these entries.

Recall the expression for the diagonal entries in~\eqref{eq:Djj_simplified}:
\begin{equation}
    D_{k,k}  = -\frac{1}{2} - \sum_{\substack{i=0 \\ i\neq k}}^n \frac{1}{x_k-x_i}.
    \label{eq:Djj_simplified2}
\end{equation}
We now simplify further by exploiting identities satisfied by the zeros of Laguerre polynomials. For Laguerre--Gauss nodes, the classical relations by~\citet{stieltjes1885} and their generalisation by~\citet{steinerberger2018} show that
\begin{equation}
    \label{eq:stieltjes}
    \sum_{\substack{i=0 \\ i\neq k}}^n \frac{1}{x_k-x_i}
    = \frac{1}{2}\left(1-\frac{\alpha+1}{x_k}\right).
\end{equation}
For augmented Laguerre--Gauss and Laguerre--Gauss--Radau nodes, the included point $x_0=0$ is not a zero of $L_n\supa$ and \eqref{eq:stieltjes} must be adjusted as follows. For $k=1,\ldots,n$ we write
\begin{equation}
    \sum_{\substack{i=0 \\ i\neq k}}^n \frac{1}{x_k-x_i}
    = \frac{1}{x_k}
      + \sum_{\substack{i=1 \\ i\neq j}}^n \frac{1}{x_k-x_i}
    = \frac{x_k+1-\alpha}{2x_k}\,, 
    \label{eq:phi_LGR}
\end{equation}
while for the augmented node itself, we utilise the identity~\citep{stieltjes1885}:
\begin{equation}
    \sum_{\substack{i=1 \\ i\neq k}}^n \frac{1}{x_0-x_i}
    = -\sum_{i=1}^n \frac{1}{x_i}
    = -\frac{n}{\alpha+1}.
    \label{eq:phi_LGR0}
\end{equation}
These expressions are substituted into~\eqref{eq:Djj_simplified2} to obtain simple, direct formulas for each of the distributions. 

\subsection{Full formulation for first-order differentiation matrices} 

The resulting compact expressions for the first-order Laguerre pseudospectral differentiation matrices of size $(n+1)\times(n+1)$ are summarised in Table~\ref{tab:difmat} for both the standard and augmented (including Radau) distributions. For the standard distribution, the indices $k,j = 1,\ldots,n+1$ are used since the collocation points $\{x_i\}_{i=1}^{n+1}$ are the roots of $L_{n+1}\supa$. For the augmented distributions, we instead use $k,j = 0,\ldots,n$, where $x_0 = 0$ and $\{x_i\}_{i=1}^{n}$ are the roots of $L_{n}\supa$.

\begin{table}[h]
\centering
\renewcommand{\arraystretch}{1.4}
\setlength{\tabcolsep}{12pt}
\caption{Entries in the first-order standard and augmented Laguerre differentiation matrices}
\label{tab:difmat}
\begin{tabular}{cccccc}
\hline
\\[-3ex]
 & ${D_{k,j}}$ & ${D_{k,k}}$ & ${D_{0,0}}$ & ${\tilde{c}_j}$ & ${\tilde{c}_0}$  \\[3pt]
 
\hline
\\[-3ex]
\textbf{Standard} 
& \multirow{2}{*}{\raisebox{-4.5ex}{$\dfrac{\tilde{c}_k/\tilde{c}_j}{x_k-x_j}$}}
& $\dfrac{-1-\alpha}{2x_k}$
& N/A 
& $\widehat{L}_n\supap(x_j)$
& N/A \\[-12pt]

\footnotesize $k,j = 1,\ldots,n+1$ & & & & & \\[6pt]

\textbf{Augmented}
& 
& $\dfrac{1-\alpha}{2x_k}$
& $-\dfrac{1}{2}-\dfrac{n}{\alpha+1}$
& $x_j\,\widehat{L}_n\supap(x_j)$ 
& $\dbinom{n+\alpha}{n}$ \\[-12pt]

\footnotesize $k,j = 0,\ldots,n$ & & & & & \\[1pt]

\hline
\\[-4ex]
\end{tabular}
\end{table}

 In the augmented case, the inclusion of the endpoint $x_0=0$ leads to a different value for the first entry $D_{0,0}$, as well as adjusted expressions for the first row $D_{0,j}$ and column $D_{k,0}$ through $\tilde{c}_0$, which is defined separately. In practice, the collocation points $x_j$ and derivative values $\widehat{L}_n\supap(x_j)$ required for the evaluation of the entries in Table~\ref{tab:difmat} are computed in a numerically stable manner using the Glaser--Liu--Rokhlin algorithm (see Section~\ref{subsec:glr}).

A perhaps more subtle advantage of the direct formula for diagonal entries is that each entry $D_{k,k}$ can be evaluated using only the corresponding collocation node $x_k$. While this property was already present for the off-diagonal entries, expressions such as~\eqref{eq:Djj_simplified} require knowledge of all nodes to compute any single diagonal term. This distinction becomes particularly important in \emph{truncated} spectral collocation, as recently employed in~\citep{hale2025}, where only a subset of the original nodes is retained. In such settings, only selected entries of the differentiation matrix are needed, and the direct formulation allows these to be computed independently, without reference to all of the nodes.

\subsection{Second-order differentiation matrices}

For second-order differentiation matrices, it follows from~\eqref{eq:derdifmat} that the diagonal entries are given by
\begin{equation}
    D^{(2)}_{k,k} = \frac{1}{4} - \phi'_k(x_k) + \phi''_k(x_k)\,, 
\end{equation}
where the second derivative $\phi''_k(x_k)$ can be derived analogously to the first derivative in Section~\ref{subsec:diagonals}. For the off-diagonal entries, the recursive procedure of~\citet{welfert1997} remains the preferred approach; however, it can now be initialised using a more accurate first-order differentiation matrix. The resulting entries of the second-order Laguerre pseudospectral differentiation matrices of size $(n+1)\times(n+1)$ are summarised in Table~\ref{tab:difmat2} for both the standard and augmented node distributions. The coefficients $\tilde c_j$ and first-order matrices $D^{\scriptscriptstyle (1) \!}_{k,k}$ and $D^{ \scriptscriptstyle (1) \!}_{k,j}$ can be found in Table~\ref{tab:difmat}. 

\begin{table}[h]
\centering
\renewcommand{\arraystretch}{1.4}
\setlength{\tabcolsep}{2pt}
\caption{Entries in the second-order standard and augmented Laguerre differentiation matrices}
\label{tab:difmat2}
\begin{tabular}{ccccc}
\hline
\\[-3ex]
& ${D^{(2)}_{k,j}}$
& ${D^{(2)}_{k,k}}$ & ${D^{(2)}_{0,0}}$ & $b$ \\[3pt]
\hline
\\[-3ex]
\textbf{Standard} 
& \multirow{2}{*}{\raisebox{-5.5ex}{$\dfrac{2}{(x_k - x_j)}\left( \dfrac{\tilde c_k }{\tilde c_j}D^{(1)}_{k,k} - D^{(1)}_{k,j} \right)$}}
& \multirow{2}{*}{\raisebox{-5.5ex}{$\dfrac{1}{12} - \dfrac{2(2n+\alpha+1)\,x_k - b}{12x_k^2}$}}
& N/A
& $ 4(\alpha+1)(\alpha+2)$\\[-9pt] 

\footnotesize $k,j = 1,\ldots,n+1$ & & \\[6pt]

\textbf{Augmented}
& 
& 
& $\dfrac{1}{4} + \dfrac{n(n+\alpha+1)}{(\alpha+1)(\alpha+2)}$
& $4(\alpha+1)(\alpha-1)$ \\[-9pt]

\footnotesize $k,j = 0,\ldots,n$ & & \\[1pt]

\hline
\end{tabular}
\end{table}

\subsection{Algorithmic improvement: the modified recurrence}
\label{subsec:modified}

Recently, \citet{huang2024} proposed a modification to the standard three-term recurrence that significantly reduces round-off errors in the computation of Laguerre polynomials and Laguerre functions. In the standard recurrence relation in~\eqref{eq:rr}, the coefficient $2n+\alpha+1-x$ appears explicitly. When $x$ is small but the constant $2n+\alpha+1$ is comparatively large, the subtraction involved in forming this coefficient leads to a loss of significant digits in floating-point arithmetic, thereby degrading numerical accuracy. To mitigate this cancellation error, \citet{huang2024} introduce the difference quantity
\begin{equation}
\delta L_n\supa = L_n\supa - L_{\scriptscriptstyle n-1}^{\scriptscriptstyle \! (\alpha)}\,,  \qquad n \geq 1\,, 
\end{equation}
and derive an alternative recurrence relation formulated in terms of $\delta L_n\supa$. Their numerical experiments demonstrate that for arguments close to zero the accuracy of computed Laguerre polynomials is improved by approximately two to four decimal digits when compared with the standard recurrence.

Since the initialisation phase of the Glaser--Liu--Rokhlin algorithm relies on evaluating Laguerre functions (and their derivatives) at arguments near the origin using a recurrence relation, it is sensible to adapt this stage of the algorithm to employ the modified recurrence of~\citet{huang2024}. Doing so improves the reliability of the starting values and propagates improved accuracy to subsequent roots generated by the algorithm, as well as the derivatives evaluated at those roots.

\section{Numerical experiments}
\label{sec:experiments}

Unless otherwise stated, all reported errors are measured relative to higher-precision reference values. These reference values were computed by implementing the procedures for computing roots and differentiation matrices using 24-digit arithmetic and then converting to double precision yielding reference values accurate to at least 16 digits. All MATLAB implementations used to produce the results and figures in this work are publicly available at~\citep{laguerre_difmat}.

\subsection{Computation of roots}

\begin{figure}[h]
    \centering
    \includegraphics[height=6cm]{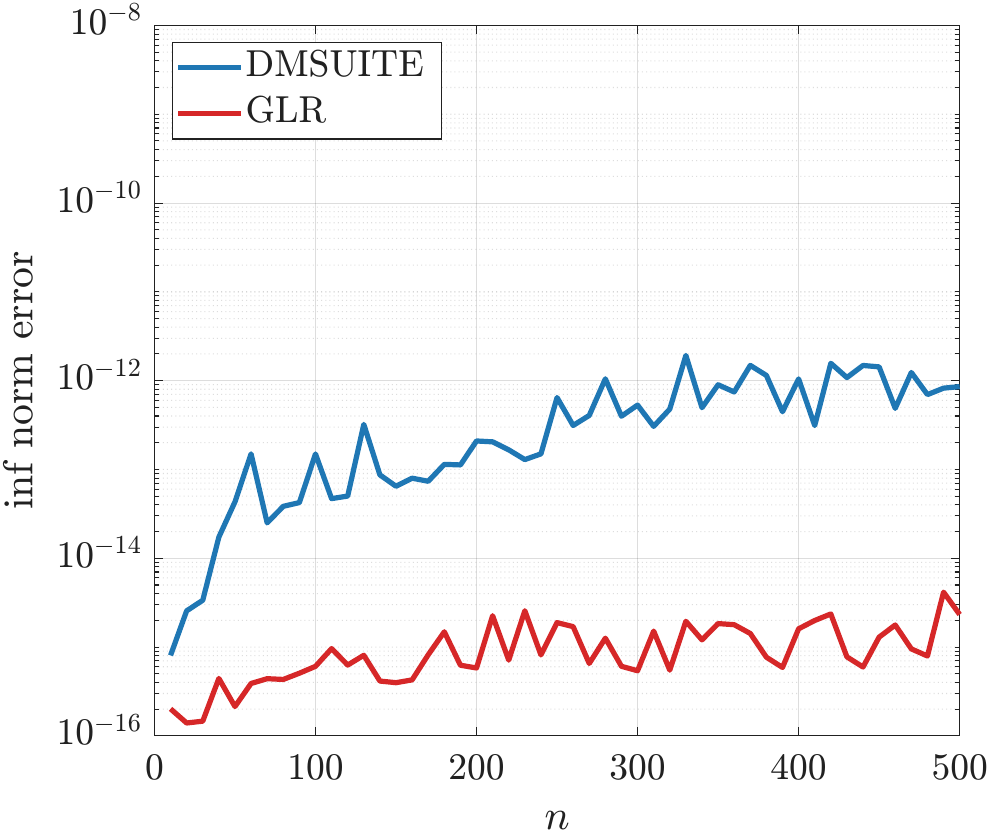}
    \hfill
    \includegraphics[height=6cm]{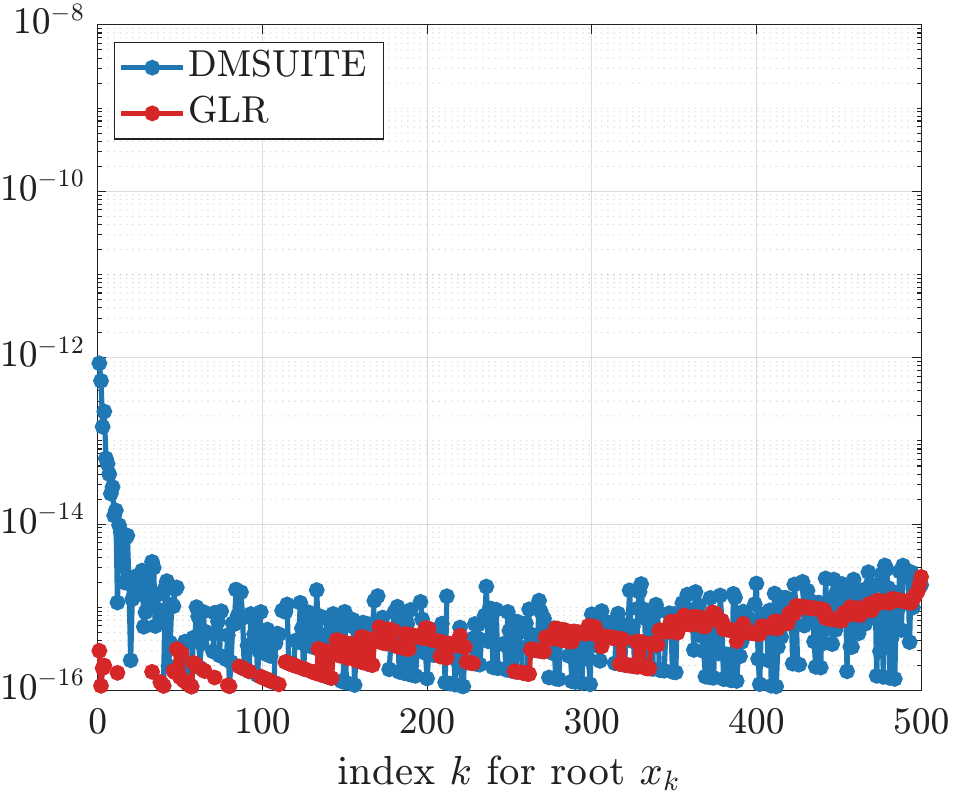}
    \caption{Relative error in computing the roots of $L_n(x)$.  Left: Infinity norm of the relative error in the computed roots as $n$ increases. The roots obtained with the Glaser--Liu--Rokhlin (GLR) algorithm exhibit smaller errors than those computed with the \texttt{DMSUITE} routine \texttt{lagroots}, which determines the roots via the Golub--Welsch algorithm.  Right: Relative error in each root $x_k$ of $L_{500}$. The Glaser--Liu--Rokhlin algorithm utilises the modified recurrence relation (Section~\ref{subsec:modified}), resulting in significantly smaller errors for small $x_k$ than those produced by \texttt{DMSUITE}.}
    \label{fig:roots_error}
\end{figure}

Figure~\ref{fig:roots_error} illustrates the performance of the Glaser--Liu--Rokhlin algorithm as a root-finding method for the classical Laguerre polynomials. The relative error in the computed roots is compared with that obtained using \texttt{DMSUITE}'s built-in root finder, \texttt{lagroots}, which determines roots via the Golub--Welsch algorithm. The left panel depicts the infinity norm of the relative error in the computed roots as the polynomial degree $n$ increases, showing that the Glaser--Liu--Rokhlin algorithm consistently achieves smaller errors than \texttt{lagroots}. The right panel shows the relative error in each root $x_k$ of $L_{500}$. In our implementation,\footnote{Our implementation is an adapted version of the \texttt{lagpts} routine from Chebfun~\citep{chebfun}.} the Glaser--Liu--Rokhlin algorithm employs the modified recurrence relation in Section~\ref{subsec:modified} rather than the standard recurrence. This modification leads to reduced error for the smallest roots $x_k$ compared with those produced by \texttt{DMSUITE}.  

\subsection{Accuracy in first-order matrix}

\begin{figure}[h!]
    \centering
    \includegraphics[height=6cm]{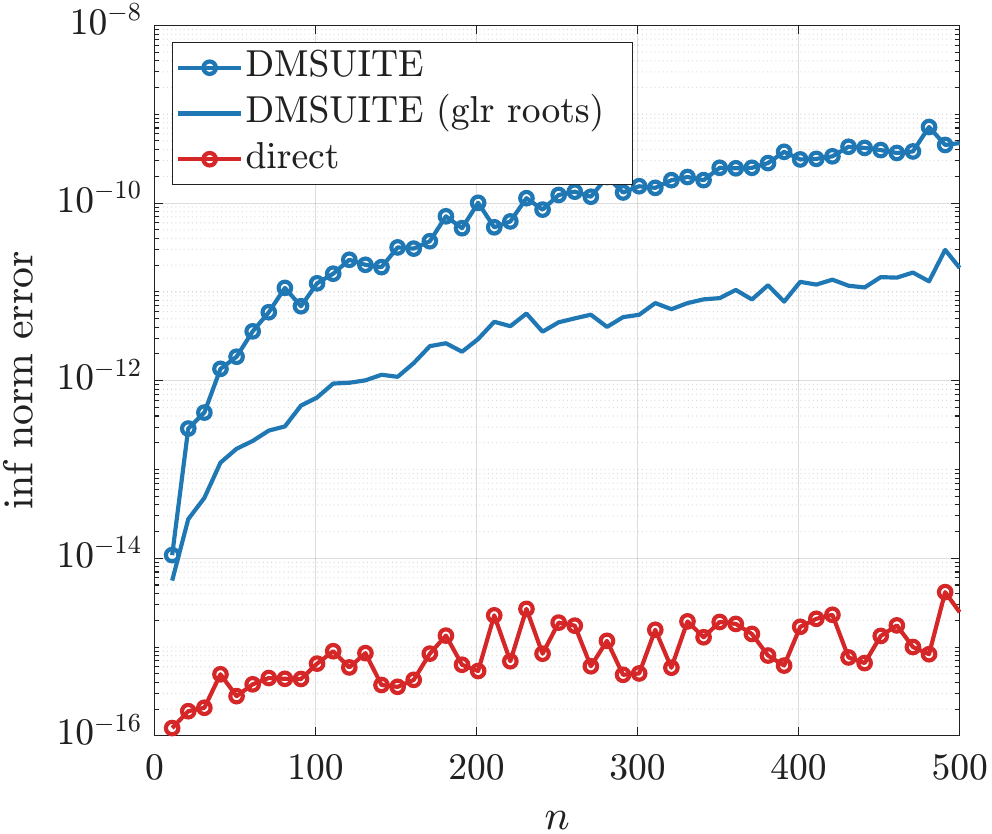}
    \hfill
    \includegraphics[height=6cm]{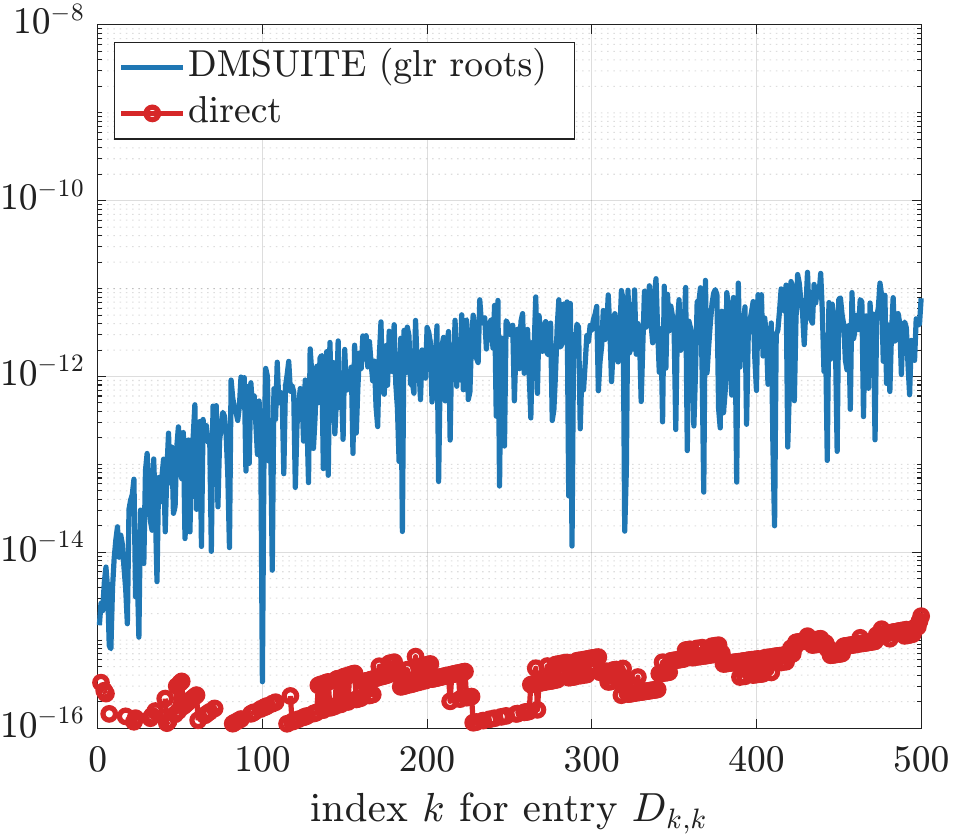}
    \caption{Relative error in the diagonal entries of the first-order augmented Laguerre--Gauss differentiation matrix. Left: Infinity norm of the relative error in the diagonal entries as $n$ increases, comparing the direct formula from Table~\ref{tab:difmat} to \texttt{DMSUITE} using its standard roots (from \texttt{lagroots}) and when more accurate roots from the Glaser--Liu--Rokhlin algorithm are used (`glr roots'). Right: Relative error in each diagonal entry for $n=500$. The direct formula produces errors comparable to the root accuracy, whereas Eq.~\eqref{eq:Djj_simplified} used by \texttt{DMSUITE} produces less accurate entries even while using the same roots as the direct formula.}
    \label{fig:error_diag_D1}
\end{figure}

The left panel of Figure~\ref{fig:breakdown} shows the relative error in the diagonal entries of the first-order augmented Laguerre--Gauss differentiation matrix. The proposed method uses the direct formula given in Table~\ref{tab:difmat} which yields smaller errors than the formula in Eq.~\eqref{eq:Djj_simplified} that is used by the \texttt{DMSUITE} implementation. Since the diagonal entries depend strongly on the roots, inaccuracies in the roots returned by \texttt{lagroots} (see Figure~\ref{fig:roots_error}) propagate into the differentiation matrix and contribute to the observed error. To separate this effect from the algorithm, we recompute the \texttt{DMSUITE} error using the more accurate roots obtained from the Glaser--Liu--Rokhlin algorithm (`glr roots'). The resulting errors are shown in Figure~\ref{fig:error_diag_D1}. While this modification improves the accuracy of the \texttt{DMSUITE} matrix entries, the error remains noticeably larger than that obtained with the direct formula. The right panel of Figure~\ref{fig:error_diag_D1} shows the relative error in each diagonal entry for $n=500$. For the direct formula, these errors are at the same level as the error in the roots shown in Figure~\ref{fig:roots_error}, whereas the \texttt{DMSUITE} formula produces substantially larger errors even while supplied with more accurate roots.

The right panel of Figure~\ref{fig:breakdown} shows the relative error in the off-diagonal entries of the first-order augmented Laguerre--Gauss differentiation matrix. The proposed method uses the scaled coefficients $\tilde{c}_j$ in the formulation described in Table~\ref{tab:difmat}, with the collocation nodes and the required Laguerre function derivatives obtained from the Glaser--Liu--Rokhlin algorithm described in Section~\ref{subsec:glr}. In contrast, the \texttt{DMSUITE} implementation evaluates the unscaled coefficients $c_j$ using~\eqref{eq:cj_def}, where the barycentric weights are computed via~\eqref{eq:vj_prod}. This procedure leads to numerical overflow (see Figure~\ref{fig:minmax}).

\begin{figure}[h]
    \centering
    \includegraphics[height = 6cm]{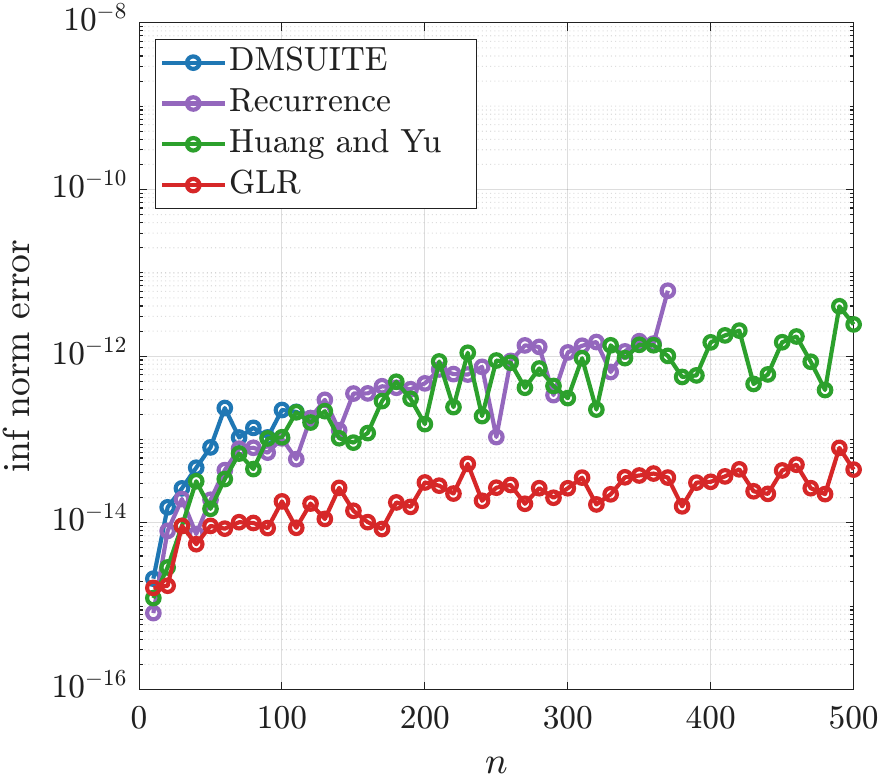}
    \hfill
    \includegraphics[height = 6cm]{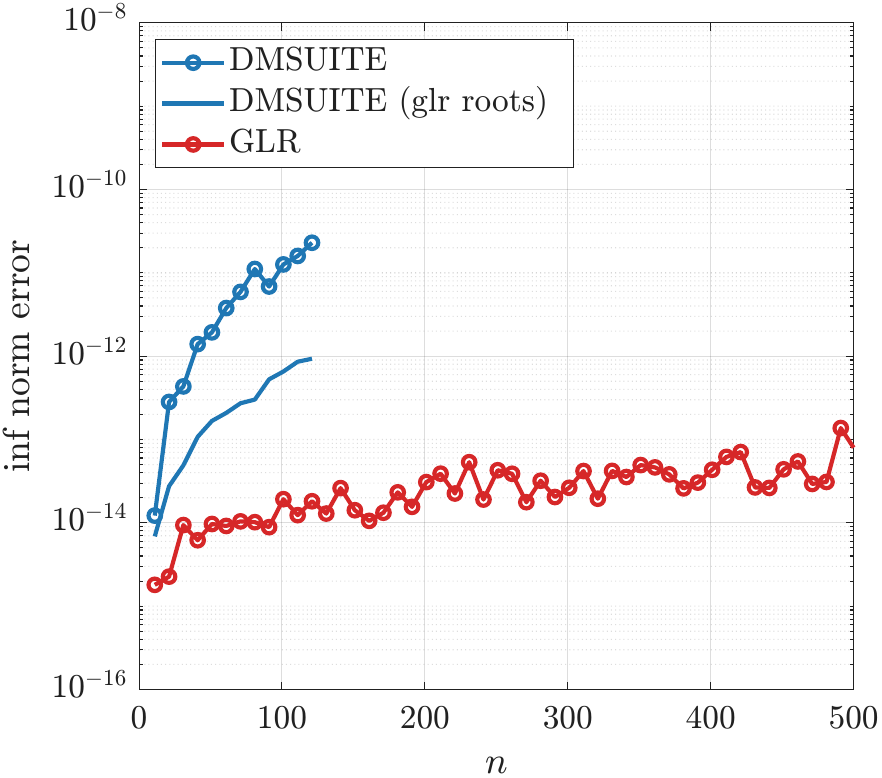}
    \caption{Left: Infinity norm of the relative error in the off-diagonal entries of the first-order augmented Laguerre--Gauss differentiation matrix. \texttt{DMSUITE} is compared to the proposed formulation for off-diagonal entries in Table~\ref{tab:difmat}, where derivatives are evaluated using three approaches: the standard Laguerre recurrence, the adaptively weighted recurrence of~\protect\citet{huang2024}, and the Glaser--Liu--Rokhlin algorithm (GLR). Although the adaptive weighting technique also remains stable for all depicted $n$, the proposed approach based on derivatives from the Glaser--Liu--Rokhlin algorithm achieves the smallest errors. Right: Infinity norm of the relative error in the off-diagonal entries of the second-order augmented Laguerre--Gauss differentiation matrix. Both methods depicted employ the recursive procedure for off-diagonal entries in Table~\ref{tab:difmat2}. The Glaser--Liu--Rokhlin based approach remains stable and achieves higher accuracy due to improved coefficient evaluation and more accurate first-order matrices. }
    \label{fig:error_offdiag_D1D2}
\end{figure}

In part, the robustness of the proposed approach stems from the derivative-based formulation in~\eqref{eq:cj_new}. Other strategies for evaluating the Laguerre function derivatives can also be used, such as the standard recurrence relation or the adaptively weighted recurrence proposed by~\citet{huang2024}. The error in computing off-diagonal entries using these alternatives are demonstrated on the left in Figure~\ref{fig:error_offdiag_D1D2}. The standard recurrence fails at $n \approx 383$ due to underflow in the explicitly evaluated exponential factor. The adaptively weighted recurrence remains stable for large $n$, but still produces larger errors than the derivatives obtained using the Glaser--Liu--Rokhlin algorithm.

An additional advantage of our approach above these alternatives lies in its more holistic use of the Glaser--Liu--Rokhlin algorithm. The algorithm already provides an efficient and competitive procedure for computing the Laguerre collocation nodes, and in doing so simultaneously yields the derivatives required for constructing the differentiation matrix. As a result, node generation and matrix assembly are performed within a single computational procedure.

\subsection{Accuracy in second-order matrix}

\begin{figure}[h]
    \centering
    \includegraphics[height=6cm]{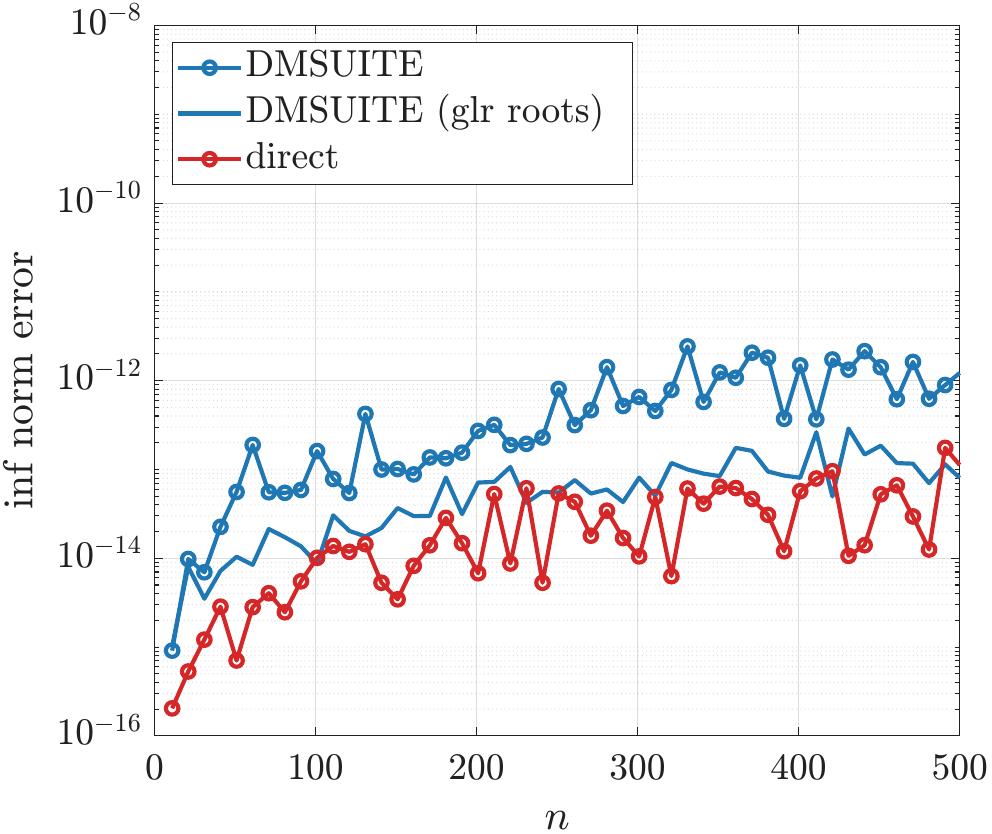}
    \hfill
    \includegraphics[height=6cm]{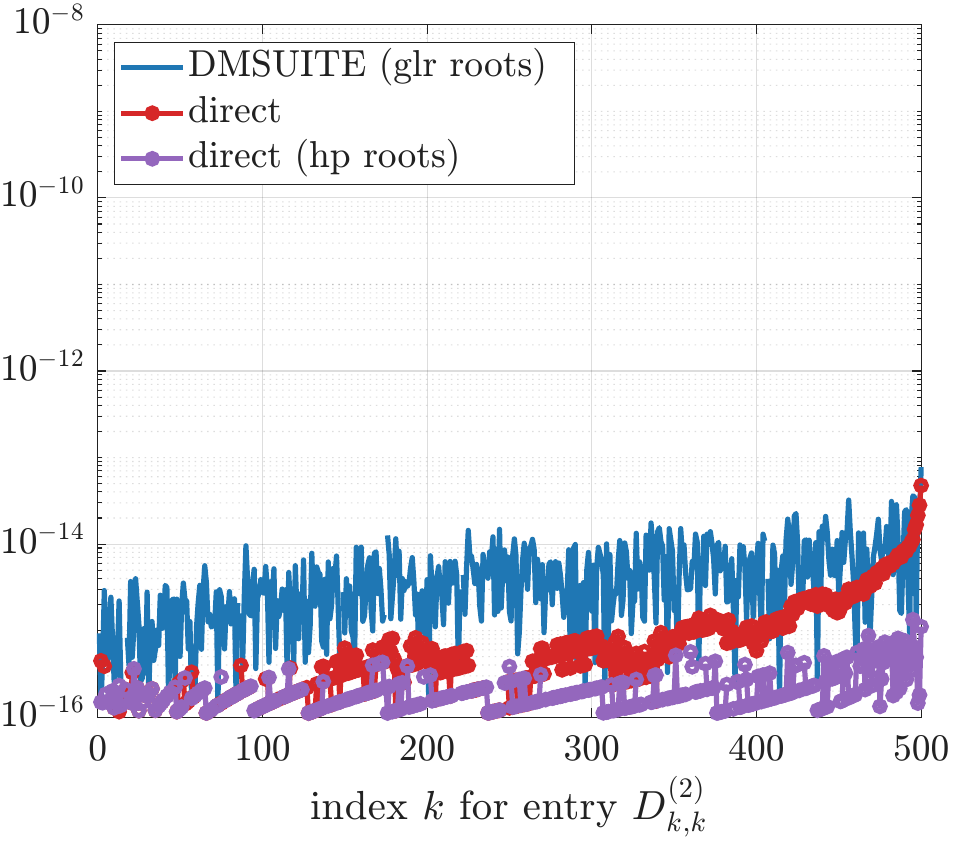}
    \caption{Relative error in the diagonal entries of the second-order augmented Laguerre--Gauss differentiation matrix. Left: Infinity norm of the relative error in the diagonal entries as $n$ increases, comparing the direct formula from Table~\ref{tab:difmat2} to \texttt{DMSUITE} using its standard roots (from \texttt{lagroots}) and when more accurate roots from the Glaser--Liu--Rokhlin algorithm are used (`glr roots'). The errors are much closer than for the first-order matrix, with the direct formula more accurate than the standard \texttt{DMSUITE} method, but only slightly better than \texttt{DMSUITE} with improved roots. Right: Relative error in each diagonal entry for $n=500$. It is clear that the highest errors in the direct formula correspond to the largest roots, where cancellation errors are amplified. The direct formula is also evaluated using high-precision roots (`hp roots'), which leads to an improvement for the largest roots. This shows the best accuracy that can be obtained using the direct formula in double precision.}
    \label{fig:error_diag_D2}
\end{figure}

Figure~\ref{fig:error_diag_D2} depicts the error in the diagonal entries of the second-order augmented Laguerre--Gauss differentiation matrix. The direct approach proposed in this work evaluates these entries using the formula given in Table~\ref{tab:difmat2}. For indices $k \geq 1$, however, the expression for $D^{\scriptscriptstyle (2)}_{\, k,k}$ is susceptible to cancellation errors at the largest roots where $x \approx 4n$. This behaviour is evident in the right panel of Figure~\ref{fig:error_diag_D2}, where the error associated with the direct method gradually increases in the tail. In the same plot we also evaluate the direct formula in double precision using high-precision Laguerre roots (`hp roots'). This leads to a modest improvement of approximately one to two orders of magnitude, which represents the practical limit of accuracy using the direct formula when computations are performed in double precision. The left and right panels also display the error obtained with \texttt{DMSUITE}, which implements the recursive scheme of~\citet{welfert1997}. As observed for the first-order differentiation matrix, the accuracy of \texttt{DMSUITE} can be improved by more accurate roots. The resulting infinity-norm errors of \texttt{DMSUITE} (left) become comparable to those of the direct method, but the relative error across most entries (right) remains approximately one order of magnitude larger.

The right panel in Figure~\ref{fig:error_offdiag_D1D2} presents the error in the off-diagonal entries of the second-order augmented Laguerre--Gauss differentiation matrix. In both approaches, the recursive procedure from Table~\ref{tab:difmat2} is used to construct the higher-order matrices. However, the Glaser--Liu--Rokhlin method yields consistently more accurate entries, owing to the improved evaluation of the coefficients $\tilde{c}_j$ and the increased accuracy of the underlying first-order differentiation matrices. As before, \texttt{DMSUITE} benefits from the use of more accurate roots. Nevertheless, the implementation breaks when the computation of the off-diagonal entries in the first-order matrix becomes unstable. 

\subsection{Example: second-order equation on the half-line}
\label{subsec:example_de}

We consider a model problem commonly used in literature to assess the performance of Laguerre spectral methods~\citep{huang2024}:
\begin{equation}
    -u''(x) + \gamma u(x) = f(x), \quad u(0) = 0, \quad \lim_{x \to \infty} u(x) = 0.
    \label{eq:ode}
\end{equation}
We prescribe the exact solution $u(x) = \sin(2x)\,e^{-x/4}$, and take $\gamma = 2$, and define the right-hand side $f(x)$ accordingly. It is well known that the accuracy of spectral methods on unbounded domains can be significantly improved  by introducing an appropriate scaling of the independent variable~\citep{shen2009,shen2000,boyd2001}. To this end, we introduce $\beta > 0$ and perform the change of variables 
\begin{equation}
    x = \frac{\tilde{x}}{\beta}\,, 
\end{equation}
where $\tilde{x}$ denotes the augmented Laguerre--Gauss nodes. Under this transformation, derivatives with respect to $\tilde{x}$ are related to derivatives with respect to $x$ by
\begin{equation}
    \frac{d}{dx} = \beta \frac{d}{d\tilde{x}}\,, 
    \qquad 
    \frac{d^2}{dx^2} = \beta^2 \frac{d^2}{d\tilde{x}^2}.
\end{equation}
In practice, this corresponds to scaling the differentiation matrices so that the first- and second-order operators are given by $\beta \, D^{(1)}$ and $\beta^2 \, D^{(2)}$, respectively. Discretising~\eqref{eq:ode} at the scaled nodes then yields the $n \times n$ linear system
\begin{equation}
    \left(-\beta^2D^{(2)} + \gamma I \right)\bm{u} = \bm{f}\,, 
\end{equation}
where $\bm{u}$ and $\bm{f}$ denote the vectors of solution and forcing values at the scaled nodes. The boundary condition at $x = 0$ is enforced by replacing the first row of the system to impose $u(0) = 0$. The decay condition at infinity is automatically satisfied by the Laguerre basis and is therefore not imposed explicitly.

\begin{figure}[h]
    \centering
    \includegraphics[height = 5.5cm]{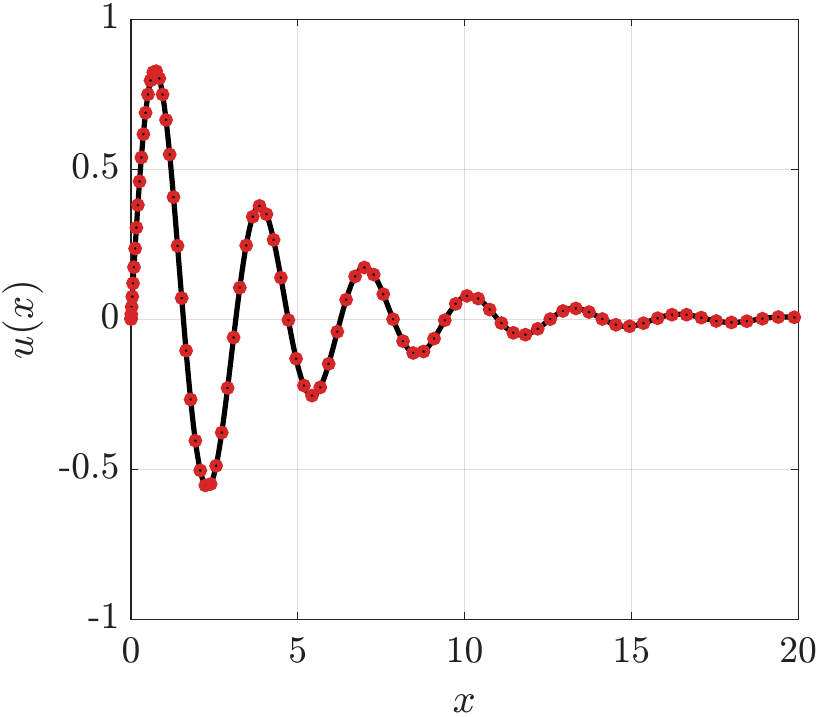}
    \hspace{5pt}
    \includegraphics[height = 5.5cm]{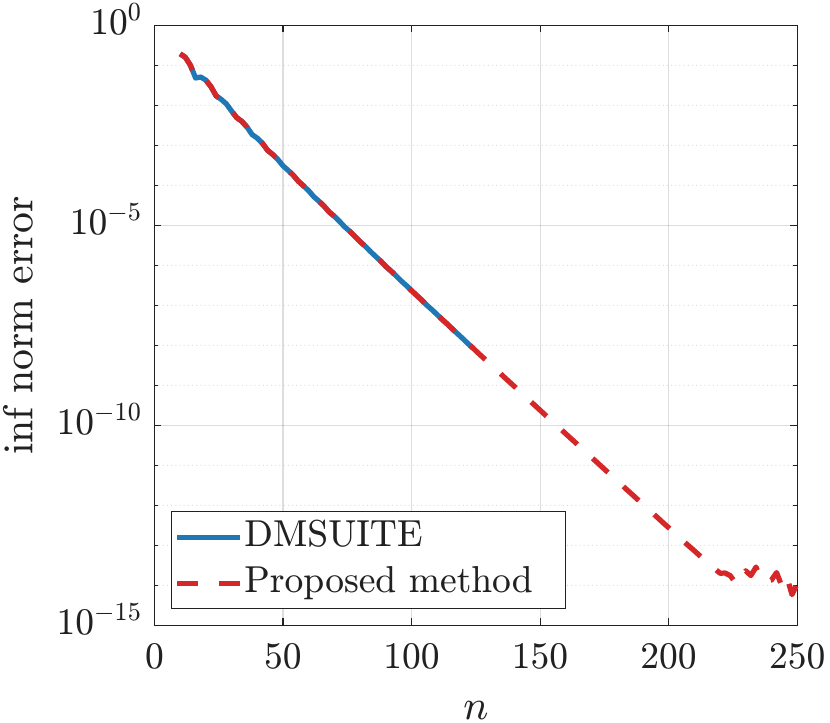}
    \caption{
    Left: Exact solution $u(x) = \sin(2x)\,e^{-x/4}$ together with the distribution of $n=230$ collocation nodes after scaling by $\beta = 4.03$. Right: Infinity norm of the absolute error of the Laguerre spectral method as the number of collocation points $n$ increases. The proposed method for generating differentiation matrices allows the solution to be resolved using $n \approx 230$ points, which would not have been possible using \texttt{DMSUITE}.
    }
    \label{fig:ode_example}
\end{figure}

The exact solution is shown on the left in Figure~\ref{fig:ode_example}, together with the error for increasing $n$ on the right. It is well known that the choice of $\beta$ has a significant influence on accuracy. For solutions of the type $e^{zx}$ where $z \in \mathbb{C}$ and $\Re(z) < 0$, an estimate for the optimal scaling is given by~\citep{huang2024}
\begin{equation}
    \beta^* = 2|z| = 2\sqrt{2^2 + (1/4)^2} \approx 4.03.
\end{equation}
Even with near-optimal scaling, a relatively large number of collocation points is required to resolve the solution, with more than $n = 200$ points needed to reach near machine precision. In this regime, it is infeasible to construct the second-order differentiation matrix using \texttt{DMSUITE}. The proposed method, in contrast, allows the solution to be resolved to errors on the order of $10^{-15}$.

\subsection{Example: eigenvalues of the time-independent Schr\"{o}dinger equation}

An important problem in quantum mechanics is the computation of eigenvalues $\lambda$ and eigenstates $y(x)$ of the time-independent Schr\"{o}dinger equation
\begin{equation}
    -y''(x) + y(x) = \lambda q(x)y(x), \quad y(0) = 0, \quad \lim_{x \to \infty} y(x) = 0.
    \label{eq:schrod}
\end{equation}
Here, the potential function $q(x)$ is taken to be the Woods--Saxon potential
\begin{equation}
    q(x) = \frac{1}{1 + e^{(x - R)/a}}\,, 
    \label{eq:potential}
\end{equation}
which is commonly used to model the interaction of a neutron with a heavy nucleus. The parameter $R$ represents the nuclear radius, $a$ characterises the surface thickness of the nucleus, and $x$ denotes the distance of the neutron from the centre.

We discretise the problem using augmented Laguerre--Gauss nodes scaled by $\beta = 10$ (as described in Section~\ref{subsec:example_de}). The differential equation is then approximated by the $n \times n$ generalised eigenvalue problem
\begin{equation}
    (-\beta^2 D^{(2)} + I)\bm{y} = \lambda Q \bm{y} \,,
\end{equation}
where $\bm{y}$ contains the approximate eigenfunction values at the nodes, $I$ is the identity matrix, and
\begin{equation}
    Q = \operatorname{diag}\left( \frac{1}{1 + e^{(x_k - R)/a}} \right).
\end{equation}
The boundary condition at $x = 0$ is imposed by forming a reduced system in which the first row and column of all the matrices are removed.

\begin{figure}[h]
    \centering
    \includegraphics[height = 5.5cm]{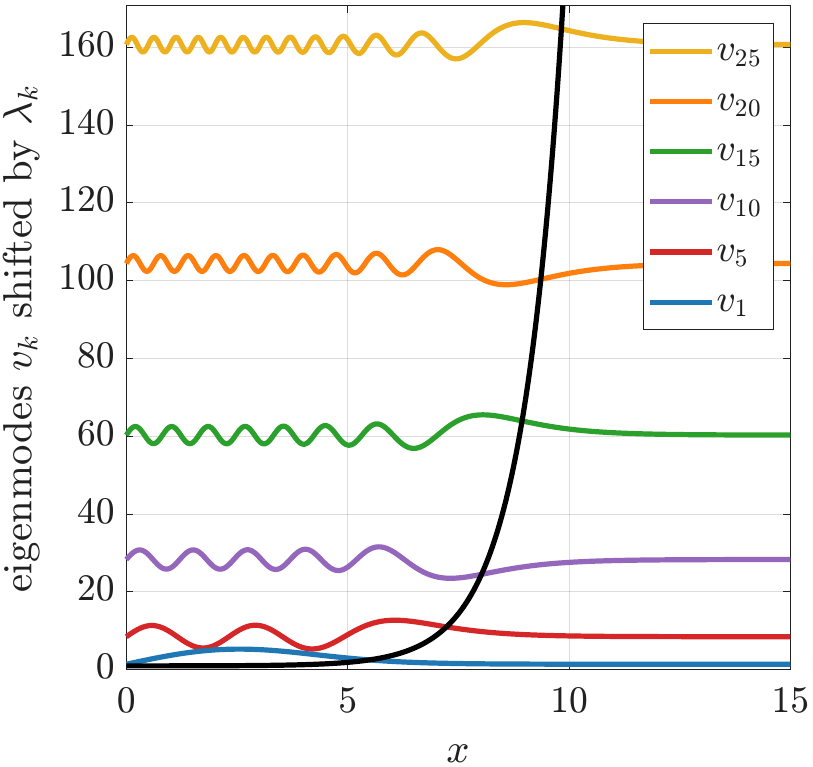}
    \hspace{10pt}
    \includegraphics[height = 5.6cm]{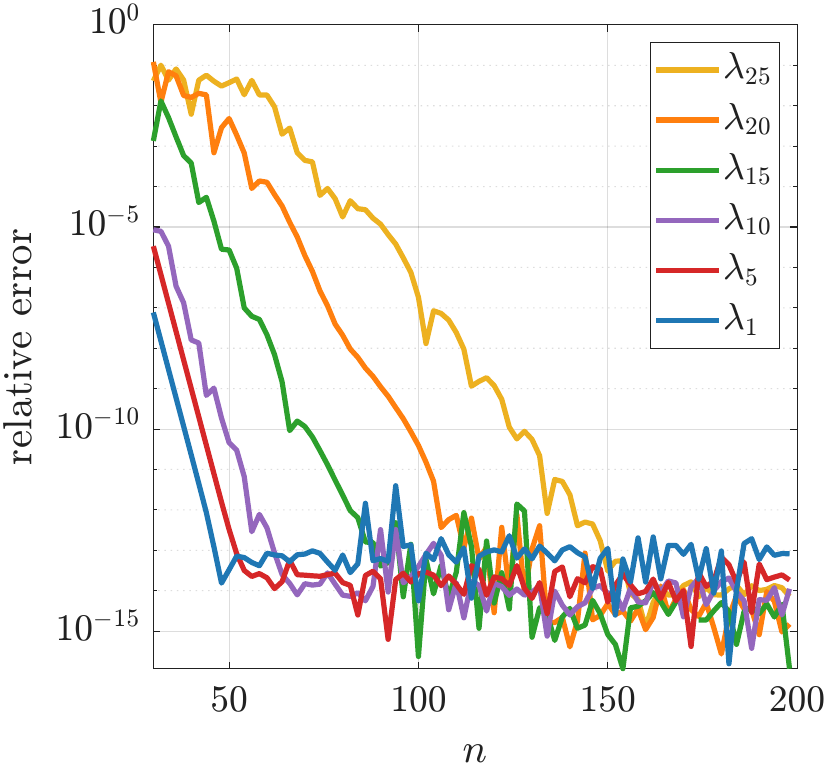}
    \caption{
    Left: Eigenfunctions of the Schr\"{o}dinger equation~\eqref{eq:schrod} computed using the Laguerre collocation method with $n = 200$ and $\beta = 10$. Each eigenmode is shifted vertically by its corresponding eigenvalue and the Woods--Saxon potential is shown in black. 
    Right: Convergence of selected eigenvalues as the number of collocation points $n$ increases, with the solution at $n = 200$ taken as the reference value. Exponential convergence is observed, with approximately 150 points needed to resolve the largest shown eigenvalue.
    }
    \label{fig:example_schrod}
\end{figure}

The eigenvalue of smallest magnitude corresponds to the ground state (often referred to as the $1s$ state) of the neutron, and is typically the physically interesting eigenvalue. The Laguerre method computes this eigenvalue to near machine precision with approximately $n \approx 50$ collocation points. To study convergence, we compute several eigenvalues for $n = 30, \dots, 200$, and in each case take the value obtained with $n = 200$ as the reference. The convergence of these eigenvalues is exponential, as shown on the right of Figure~\ref{fig:example_schrod}. To resolve the eigenvalue $\lambda_{25}$ we require more than 150 collocation points. On the left of Figure~\ref{fig:example_schrod}, we plot the corresponding eigenfunctions, each shifted vertically by an amount equal to its eigenvalue. The black curve represents the Woods--Saxon potential. As expected, the functions oscillate within the potential well and rapidly decay to zero outside of it.

\section{Conclusion}
\label{sec:conclusion}

While the differentiation matrices in Laguerre spectral collocation are typically well-scaled, classical constructions thereof are prone to severe numerical instability. This instability arises from the rapid growth of Laguerre polynomials and their derivatives, the exponential decay of the associated Laguerre functions, and the cancellation of extreme intermediate quantities. These effects lead to overflow, underflow, and loss of significance in standard implementations, consequently limiting the practical range of Laguerre collocation methods.

The key contribution of this paper is a reformulation of the differentiation matrix entries that avoids these numerical difficulties. For the off-diagonal entries, we suggest a representation in terms of derivatives of Laguerre polynomials and emphasise that all required quantities must be computed simultaneously to avoid the computation of extremely large or small values. By exploiting the fast algorithm of~\citet{glaser2007}, we obtain the collocation nodes and the associated derivative information in a unified and stable manner. For the diagonal entries, we advocate a direct closed-form expression and demonstrate that it yields improved accuracy compared to alternatives. 

Numerical experiments confirm that the proposed approach extends the range of computability beyond that of widely used classical implementations. In particular, the stable construction avoids the breakdown observed in the standard software package \texttt{DMSUITE}~\citep{dmsuite} and maintains high accuracy for substantially larger values of $n$. The resulting all-in-one procedure is fast, robust, and straightforward to implement.

Several avenues for future work remain. The ideas developed here could be extended to interpolation frameworks by incorporating the exponential factor directly into \emph{barycentric resampling matrices}, such as those of~\citet{hale2025}. Similarly, the stable evaluation of Laguerre-based quadrature rules and their associated weights warrants a comparable level of attention.

An investigation into strategies for enforcing boundary conditions would also be of interest. For example, it remains to be determined whether augmented Laguerre--Gauss nodes combined with a row-replacement strategy (as in Section~\ref{subsec:example_de}) offer advantages over standard Laguerre--Gauss nodes with a basis recombination approach. Alternative approaches could also be compared, including augmenting the system with additional constraints and solving a least-squares problem, or employing rectangular projection techniques~\citep{driscoll2016}, in which the system is projected onto a lower-degree subspace before being augmented to recover a square system.

Another important direction, which will be explored in future work, concerns improving the convergence of these methods. While scaling factors, as introduced in Section~\ref{subsec:example_de}, play a crucial role, no general guidelines currently exist for selecting optimal scaling parameters for a given problem. Developing such criteria therefore remains a key challenge. In addition, \emph{truncated} spectral collocation methods, in which only a subset of the original nodes is retained, show promise for enhancing convergence rates in certain settings. In this context, the direct formula for the diagonal entries is particularly advantageous, as it enables the efficient construction of submatrices without requiring the full differentiation matrix.

\bibliographystyle{apalike}
\bibliography{reference}

\begin{thebibliography}{}

\bibitem[Adams, 2022]{julia_dmsuite}
Adams, L. (2022).
\newblock {DMSuite.jl}.
\newblock \url{https://github.com/l90lpa/DMSuite.jl}.
\newblock GitHub repository, accessed February 3, 2026.

\bibitem[Adhikari, 2013]{python_dmsuite}
Adhikari, R. (2013).
\newblock pyddx.
\newblock \url{https://github.com/ronojoy/pyddx}.
\newblock GitHub repository, accessed February 3, 2026.

\bibitem[Baltensperger and Trummer, 2003]{baltensperger2003}
Baltensperger, R. and Trummer, M.~R. (2003).
\newblock Spectral differencing with a twist.
\newblock {\em SIAM J. Sci. Comput.}, 24(5):1465--1487.

\bibitem[Boyd, 2001]{boyd2001}
Boyd, J.~P. (2001).
\newblock {\em Chebyshev and {F}ourier spectral methods}.
\newblock Dover Publications, 2nd edition.

\bibitem[Burns et~al., 2020]{dedalus}
Burns, K.~J., Vasil, G.~M., Oishi, J.~S., Lecoanet, D., and Brown, B.~P. (2020).
\newblock Dedalus: A flexible framework for numerical simulations with spectral methods.
\newblock {\em Phys. Rev. Res.}, 2(2).

\bibitem[Canuto et~al., 1988]{canuto1988}
Canuto, C., Hussaini, M.~Y., Quarteroni, A., and Zang, T.~A. (1988).
\newblock {\em Spectral methods in fluid dynamics}.
\newblock Springer Series in Computational Physics. Springer-Verlag, New York.

\bibitem[Canuto et~al., 2006]{canuto2006}
Canuto, C., Hussaini, M.~Y., Quarteroni, A., and Zang, T.~A. (2006).
\newblock {\em Spectral methods: Fundamentals in single domains}.
\newblock Scientific Computation. Springer-Verlag, Berlin.

\bibitem[Driscoll and Hale, 2016]{driscoll2016}
Driscoll, T.~A. and Hale, N. (2016).
\newblock Rectangular spectral collocation.
\newblock {\em IMA J. Numer. Anal.}, 36(1):108--132.

\bibitem[Driscoll et~al., 2014]{chebfun}
Driscoll, T.~A., Hale, N., and Trefethen, L.~N., editors (2014).
\newblock {\em Chebfun Guide}.
\newblock Pafnuty Publications, Oxford.

\bibitem[Funaro, 1990]{funaro1990}
Funaro, D. (1990).
\newblock Computational aspects of pseudospectral {L}aguerre approximations.
\newblock {\em Appl. Numer. Math.}, 6(6):447--457.

\bibitem[Gatteschi, 2002]{gatteschi2002}
Gatteschi, L. (2002).
\newblock Asymptotics and bounds for the zeros of {L}aguerre polynomials: a survey.
\newblock {\em J. Comput. Appl. Math.}, 144(1-2):7--27.

\bibitem[Gautschi, 2004]{gautschi2004}
Gautschi, W. (2004).
\newblock {\em Orthogonal polynomials: computation and approximation}.
\newblock Numerical Mathematics and Scientific Computation. Oxford University Press, New York.

\bibitem[Gheorghiu, 2013]{gheorghiu2013}
Gheorghiu, C.-I. (2013).
\newblock {Laguerre} collocation solutions to boundary layer type problems.
\newblock {\em Numer. Algorithms}, 64(2):385--401.

\bibitem[Gil et~al., 2018]{gil2018}
Gil, A., Segura, J., and Temme, N.~M. (2018).
\newblock Asymptotic approximations to the nodes and weights of {G}auss-{H}ermite and {G}auss-{L}aguerre quadratures.
\newblock {\em Stud. Appl. Math.}, 140(3):298--332.

\bibitem[Gil et~al., 2019]{gil2019}
Gil, A., Segura, J., and Temme, N.~M. (2019).
\newblock Fast, reliable and unrestricted iterative computation of {G}auss-{H}ermite and {G}auss-{L}aguerre quadratures.
\newblock {\em Numer. Math.}, 143(3):649--682.

\bibitem[Gil et~al., 2025]{gil2025}
Gil, A., Segura, J., and Temme, N.~M. (2025).
\newblock Fast and accurate computation of classical {Gaussian} quadratures.
\newblock Preprint, arXiv:2509.16716.

\bibitem[Glaser et~al., 2007]{glaser2007}
Glaser, A., Liu, X., and Rokhlin, V. (2007).
\newblock A fast algorithm for the calculation of the roots of special functions.
\newblock {\em SIAM J. Sci. Comput.}, 29(4):1420--1438.

\bibitem[Golub and Welsch, 1969]{golub1969}
Golub, G.~H. and Welsch, J.~H. (1969).
\newblock Calculation of {G}auss quadrature rules.
\newblock {\em Math. Comp.}, 23:221--230.

\bibitem[Gottlieb and Orszag, 1977]{gottlieb1977}
Gottlieb, D. and Orszag, S.~A. (1977).
\newblock {\em Numerical analysis of spectral methods: Theory and applications}, volume~26 of {\em CBMS-NSF Regional Conference Series in Applied Mathematics}.
\newblock SIAM.

\bibitem[Gu and Eisenstat, 1995]{gu1995}
Gu, M. and Eisenstat, S.~C. (1995).
\newblock A divide-and-conquer algorithm for the symmetric tridiagonal eigenproblem.
\newblock {\em SIAM J. Matrix Anal. Appl.}, 16(1):172--191.

\bibitem[Hale et~al., 2025]{hale2025}
Hale, N., Thomann, E., and Weideman, J. (2025).
\newblock Approximate solutions to a nonlinear functional differential equation.
\newblock Preprint, arXiv:2401.11733.

\bibitem[Henrici, 1964]{henrici1964}
Henrici, P. (1964).
\newblock {\em Elements of numerical analysis}.
\newblock John Wiley \& Sons.

\bibitem[Huang and Yu, 2024]{huang2024}
Huang, S. and Yu, H. (2024).
\newblock Improved {L}aguerre spectral methods with less round-off errors and better stability.
\newblock {\em J. Sci. Comput.}, 101(3).

\bibitem[Jewell, 2009]{jewell2026}
Jewell, N. (2009).
\newblock {Laguerre} spectral/pseudospectral library.
\newblock \url{https://www.mathworks.com/matlabcentral/fileexchange/26089-laguerre-spectral-pseudospectral-library}.
\newblock MATLAB Central File Exchange, accessed February 3, 2026.

\bibitem[Latifi and Delkhosh, 2019]{spsmat}
Latifi, S. and Delkhosh, M. (2019).
\newblock {SPSMAT}: {GNU Octave} software package for spectral and pseudospectral methods.

\bibitem[Maday et~al., 1985]{maday1985}
Maday, Y., Pernaud-Thomas, B., and Vandeven, H. (1985).
\newblock Reappraisal of {Laguerre} type spectral methods.
\newblock {\em La Recherche Aerospatiale}, 6:13--35.

\bibitem[Mortensen, 2018]{shenfun}
Mortensen, M. (2018).
\newblock Shenfun: High performance spectral {Galerkin} computing platform.
\newblock {\em Journal of Open Source Software}, 3(31):1071.

\bibitem[Nel, 2026]{laguerre_difmat}
Nel, E. (2026).
\newblock {MATLAB} code to accompany: Construction of {L}aguerre pseudospectral differentiation matrices.
\newblock \url{https://github.com/Emma-Nel/LaguerreDifmat}.
\newblock GitHub repository, accessed April 21, 2026.

\bibitem[Olver et~al., 2010]{nist2010}
Olver, F. W.~J., Lozier, D.~W., Boisvert, R.~F., and Clark, C.~W., editors (2010).
\newblock {\em N{IST} handbook of mathematical functions}.
\newblock Cambridge University Press.

\bibitem[Olver and Townsend, 2014]{approxfun}
Olver, S. and Townsend, A. (2014).
\newblock A practical framework for infinite-dimensional linear algebra.
\newblock In {\em 2014 First Workshop for High Performance Technical Computing in Dynamic Languages}, pages 57--62.

\bibitem[Opsomer and Huybrechs, 2023]{opsomer2023}
Opsomer, P. and Huybrechs, D. (2023).
\newblock High-order asymptotic expansions of {G}aussian quadrature rules with classical and generalized weight functions.
\newblock {\em J. Comput. Appl. Math.}, 434.

\bibitem[Shen, 2000]{shen2000}
Shen, J. (2000).
\newblock Stable and efficient spectral methods in unbounded domains using {L}aguerre functions.
\newblock {\em SIAM J. Numer. Anal.}, 38(4):1113--1133.

\bibitem[Shen et~al., 2011]{shen2011}
Shen, J., Tang, T., and Wang, L.-L. (2011).
\newblock {\em Spectral methods: Algorithms, analysis and applications}, volume~41 of {\em Springer Series in Computational Mathematics}.
\newblock Springer, Heidelberg.

\bibitem[Shen and Wang, 2009]{shen2009}
Shen, J. and Wang, L.-L. (2009).
\newblock Some recent advances on spectral methods for unbounded domains.
\newblock {\em Commun. Comput. Phys.}, 5(2-4):195--241.

\bibitem[Steinerberger, 2018]{steinerberger2018}
Steinerberger, S. (2018).
\newblock Electrostatic interpretation of zeros of orthogonal polynomials.
\newblock {\em Proc. Amer. Math. Soc.}, 146(12):5323--5331.

\bibitem[Stieltjes, 1885]{stieltjes1885}
Stieltjes, T.~J. (1885).
\newblock Sur certains polyn{\^o}mes qui v{\'e}rifient une {\'e}quation diff{\'e}rentielle lin{\'e}aire du second ordre et sur la th{\'e}orie des fonctions de lam{\'e}.
\newblock {\em Acta Mathematica}, 6:321--326.

\bibitem[Szeg\"{o}, 1939]{szego1967}
Szeg\"{o}, G. (1939).
\newblock {\em Orthogonal polynomials}, volume~23 of {\em American Mathematical Society Colloquium Publications}.
\newblock American Mathematical Society, 3rd edition.

\bibitem[von Winckel, 2004]{vonwinckel2026}
von Winckel, G. (2004).
\newblock Pseudospectral differentiation on an arbitrary grid.
\newblock \url{https://www.mathworks.com/matlabcentral/fileexchange/5515-pseudospectral-differentiation-on-an-arbitrary-grid}.
\newblock MATLAB Central File Exchange, accessed February 3, 2026.

\bibitem[Wang, 2024]{wang2024}
Wang, H. (2024).
\newblock Convergence analysis of {L}aguerre approximations for analytic functions.
\newblock {\em Math. Comp.}, 93(350):2861--2884.

\bibitem[Wang and Guo, 2008]{wang2008}
Wang, T.-J. and Guo, B.-Y. (2008).
\newblock Composite generalized {L}aguerre-{L}egendre pseudospectral method for {F}okker-{P}lanck equation in an infinite channel.
\newblock {\em Appl. Numer. Math.}, 58(10):1448--1466.

\bibitem[Wang and Sun, 2016]{wang2016}
Wang, T.-J. and Sun, T. (2016).
\newblock Mixed pseudospectral method for heat transfer.
\newblock {\em Math. Model. Anal.}, 21(2):199--219.

\bibitem[Weideman, 2003a]{weideman}
Weideman, J. A.~C. (2003a).
\newblock {A MATLAB Differentiation Matrix Suite}.
\newblock \url{https://appliedmaths.sun.ac.za/~weideman/research/differ.html}.
\newblock Personal website of J.A.C. Weideman, accessed March 4, 2026.

\bibitem[Weideman, 2003b]{dmsuite}
Weideman, J. A.~C. (2003b).
\newblock {DMSUITE}.
\newblock \url{https://www.mathworks.com/matlabcentral/fileexchange/29-dmsuite}.
\newblock MATLAB Central File Exchange, accessed March 4, 2026.

\bibitem[Weideman and Reddy, 2000]{weideman2000}
Weideman, J. A.~C. and Reddy, S.~C. (2000).
\newblock A {MATLAB} differentiation matrix suite.
\newblock {\em ACM Trans. Math. Software}, 26(4):465--519.

\bibitem[Welfert, 1997]{welfert1997}
Welfert, B.~D. (1997).
\newblock Generation of pseudospectral differentiation matrices {I}.
\newblock {\em SIAM J. Numer. Anal.}, 34(4):1640--1657.

\end{thebibliography}

\end{document}